\documentclass[12pt]{article}
\usepackage{amsmath,amssymb,amsthm,amscd,a4wide}

\begin{document}

\newcommand{\End}{{\rm{End}\ts}}
\newcommand{\Hom}{{\rm{Hom}}}
\newcommand{\ch}{{\rm{ch}\ts}}
\newcommand{\non}{\nonumber}
\newcommand{\wt}{\widetilde}
\newcommand{\wh}{\widehat}
\newcommand{\ot}{\otimes}
\newcommand{\la}{\lambda}
\newcommand{\La}{\Lambda}
\newcommand{\De}{\Delta}
\newcommand{\al}{\alpha}
\newcommand{\be}{\beta}
\newcommand{\ga}{\gamma}
\newcommand{\ve}{\varepsilon^{}}
\newcommand{\ep}{\epsilon}
\newcommand{\ka}{\kappa}
\newcommand{\vk}{\varkappa}
\newcommand{\si}{\sigma}
\newcommand{\vp}{\varphi}
\newcommand{\de}{\delta^{}}
\newcommand{\ze}{\zeta}
\newcommand{\om}{\omega}
\newcommand{\hra}{\hookrightarrow}
\newcommand{\ee}{\ep^{}}
\newcommand{\su}{s^{}}
\newcommand{\ts}{\,}
\newcommand{\vac}{\mathbf{1}}
\newcommand{\di}{\partial}
\newcommand{\qin}{q^{-1}}
\newcommand{\tss}{\hspace{1pt}}
\newcommand{\Sr}{ {\rm S}}
\newcommand{\U}{ {\rm U}}
\newcommand{\X}{ {\rm X}}
\newcommand{\Y}{ {\rm Y}}
\newcommand{\Z}{{\rm Z}}
\newcommand{\AAb}{\mathbb{A}\tss}
\newcommand{\CC}{\mathbb{C}\tss}
\newcommand{\QQ}{\mathbb{Q}\tss}
\newcommand{\SSb}{\mathbb{S}\tss}
\newcommand{\ZZ}{\mathbb{Z}\tss}
\newcommand{\Ac}{\mathcal{A}}
\newcommand{\Lc}{\mathcal{L}}
\newcommand{\Mc}{\mathcal{M}}
\newcommand{\Pc}{\mathcal{P}}
\newcommand{\Qc}{\mathcal{Q}}
\newcommand{\Tc}{\mathcal{T}}
\newcommand{\Sc}{\mathcal{S}}
\newcommand{\Bc}{\mathcal{B}}
\newcommand{\Dc}{\mathcal{D}}
\newcommand{\Ec}{\mathcal{E}}
\newcommand{\Fc}{\mathcal{F}}
\newcommand{\Hc}{\mathcal{H}}
\newcommand{\Uc}{\mathcal{U}}
\newcommand{\Wc}{\mathcal{W}}
\newcommand{\Zc}{\mathcal{Z}}
\newcommand{\Ar}{{\rm A}}
\newcommand{\Ir}{{\rm I}}
\newcommand{\Gr}{{\rm G}}
\newcommand{\Or}{{\rm O}}
\newcommand{\Spr}{{\rm Sp}}
\newcommand{\Wr}{{\rm W}}
\newcommand{\Zr}{{\rm Z}}
\newcommand{\gl}{\mathfrak{gl}}
\newcommand{\Pf}{{\rm Pf}\ts}
\newcommand{\oa}{\mathfrak{o}}
\newcommand{\spa}{\mathfrak{sp}}
\newcommand{\g}{\mathfrak{g}}
\newcommand{\h}{\mathfrak h}
\newcommand{\n}{\mathfrak n}
\newcommand{\z}{\mathfrak{z}}
\newcommand{\Zgot}{\mathfrak{Z}}
\newcommand{\p}{\mathfrak{p}}
\newcommand{\sll}{\mathfrak{sl}}
\newcommand{\agot}{\mathfrak{a}}
\newcommand{\qdet}{ {\rm qdet}\ts}
\newcommand{\Ber}{ {\rm Ber}\ts}
\newcommand{\Mat}{{\rm{Mat}}}
\newcommand{\HC}{ {\mathcal HC}}
\newcommand{\ev}{ {\rm ev}}
\newcommand{\cdet}{ {\rm cdet}}
\newcommand{\tr}{ {\rm tr}}
\newcommand{\gr}{{\rm gr}\ts}
\newcommand{\str}{ {\rm str}}
\newcommand{\loc}{{\rm loc}}
\newcommand{\sgn}{ {\rm sgn}\ts}
\newcommand{\hb}{\mathbf{h}}
\newcommand{\ba}{\bar{a}}
\newcommand{\bb}{\bar{b}}
\newcommand{\bi}{\bar{\imath}}
\newcommand{\bj}{\bar{\jmath}}
\newcommand{\bk}{\bar{k}}
\newcommand{\bl}{\bar{l}}
\newcommand{\Sym}{\mathfrak S}
\newcommand{\fand}{\quad\text{and}\quad}
\newcommand{\Fand}{\qquad\text{and}\qquad}
\newcommand{\For}{\qquad\text{or}\qquad}
\newcommand{\vt}{{\tss|\hspace{-1.5pt}|\tss}}

\renewcommand{\theequation}{\arabic{section}.\arabic{equation}}

\newtheorem{thm}{Theorem}[section]
\newtheorem{lem}[thm]{Lemma}
\newtheorem{prop}[thm]{Proposition}
\newtheorem{cor}[thm]{Corollary}
\newtheorem{conj}[thm]{Conjecture}
\newtheorem*{mthm}{Main Theorem}

\theoremstyle{definition}
\newtheorem{defin}[thm]{Definition}

\theoremstyle{remark}
\newtheorem{remark}[thm]{Remark}
\newtheorem{example}[thm]{Example}

\newcommand{\bth}{\begin{thm}}
\renewcommand{\eth}{\end{thm}}
\newcommand{\bpr}{\begin{prop}}
\newcommand{\epr}{\end{prop}}
\newcommand{\ble}{\begin{lem}}
\newcommand{\ele}{\end{lem}}
\newcommand{\bco}{\begin{cor}}
\newcommand{\eco}{\end{cor}}
\newcommand{\bde}{\begin{defin}}
\newcommand{\ede}{\end{defin}}
\newcommand{\bex}{\begin{example}}
\newcommand{\eex}{\end{example}}
\newcommand{\bre}{\begin{remark}}
\newcommand{\ere}{\end{remark}}
\newcommand{\bcj}{\begin{conj}}
\newcommand{\ecj}{\end{conj}}

\newcommand{\bal}{\begin{aligned}}
\newcommand{\eal}{\end{aligned}}
\newcommand{\beq}{\begin{equation}}
\newcommand{\eeq}{\end{equation}}
\newcommand{\ben}{\begin{equation*}}
\newcommand{\een}{\end{equation*}}

\newcommand{\bpf}{\begin{proof}}
\newcommand{\epf}{\end{proof}}

\def\beql#1{\begin{equation}\label{#1}}

\title{\Large\bf Feigin--Frenkel center
in types $B$, $C$ and $D$}

\author{A. I. Molev}

\date{} 
\maketitle

\vspace{25 mm}

\begin{abstract}
For each simple Lie algebra $\g$ consider the corresponding affine
vertex algebra $V_{\rm{crit}}(\g)$ at the critical level.
The center of this vertex algebra is a commutative associative
algebra whose structure was described by a remarkable theorem of Feigin
and Frenkel about two decades ago. However, only recently
simple formulas for the generators of the center
were found for the Lie algebras of type $A$
following Talalaev's discovery of explicit higher Gaudin
Hamiltonians. We give explicit formulas for generators
of the centers of the affine vertex algebras $V_{\rm{crit}}(\g)$
associated with the simple Lie algebras $\g$
of types $B$, $C$ and $D$.
The construction relies on the Schur--Weyl
duality involving the Brauer algebra, and the generators
are expressed as weighted traces over tensor spaces
and, equivalently, as traces
over the spaces of singular vectors for the action
of the Lie algebra $\sll_2$ in the context of Howe duality.
This leads to explicit constructions
of commutative subalgebras of the universal enveloping
algebras $\U(\g[t])$ and $\U(\g)$, and to higher order Hamiltonians
in the Gaudin model associated with each Lie algebra $\g$.
We also introduce analogues of the Bethe
subalgebras of the Yangians $\Y(\g)$ and show
that their graded images coincide with the respective
commutative subalgebras of $\U(\g[t])$.
\end{abstract}

\vspace{45 mm}

\noindent
School of Mathematics and Statistics\newline
University of Sydney,
NSW 2006, Australia\newline
alexander.molev@sydney.edu.au

\newpage

%

\section{Introduction}
\label{sec:int}
\setcounter{equation}{0}

For each simple Lie algebra $\g$ over $\CC$ the
vacuum module $V_{\ka}(\g)$ at the level $\ka\in\CC$
over the affine Kac--Moody algebra $\wh\g$
has a vertex algebra structure. The family of
affine vertex algebras $V_{\ka}(\g)$ has profound
connections in geometry and mathematical physics;
see e.g. \cite{dk:fa}, \cite{fb:va} and \cite{k:va}.

The critical value $\ka=-h^{\vee}$, where
$h^{\vee}$ is the dual Coxeter number for $\g$, plays
a particular role, as the affine vertex algebra
at the critical level
contains a big center
$\z(\wh\g)$. The structure of the center was described by
a remarkable theorem of Feigin
and Frenkel in \cite{ff:ak} which states that $\z(\wh\g)$
is the algebra of polynomials in infinitely many variables
associated with the algebra of invariants $\Sr(\g)^{\g}$
in the symmetric algebra.
The vertex algebra structure on the vacuum module
brings up a few bridges connecting the Feigin--Frenkel
center $\z(\wh\g)$ with the representation theory of
the affine Kac--Moody algebra $\wh\g$. In particular,
the higher Sugawara operators associated with the
elements of the center provide
singular vectors in the Verma modules $M(\la)$ at
the critical level which leads to a proof of the Kac--Kazhdan
conjecture for the characters of the irreducible quotients
$L(\la)$ (for the classical types the proofs were given
in \cite{gw:ho} and \cite{h:so}). Moreover, the Sugawara
operators act as scalars in the Wakimoto modules over $\wh\g$
thus providing an affine analogue of
the Harish-Chandra homomorphism and leading to a description
of the associated (Langlands dual)
classical $\Wc$-algebra; see \cite{f:lc}
for a detailed exposition of these results,
and also \cite{af:rv}
and \cite{fg:lg} for the role of the center in
the representation theory of $\wh\g$.
A striking
connection of $\z(\wh\g)$ and the Hamiltonians
of the quantum Gaudin model was discovered in \cite{ffr:gm};
see also more recent work \cite{fft:gm} and \cite{r:si} for extensions
to some generalized Gaudin models and for related constructions of
commutative subalgebras in enveloping algebras.
This connection, together
with simple explicit formulas
of the higher Gaudin Hamiltonians by Talalaev~\cite{t:qg}
in type $A$, has produced equally simple formulas
for generators of the Feigin--Frenkel center; see
\cite{cm:ho}, \cite{ct:qs}.
A super-version of these results was given in \cite{mr:mm}
together with simpler arguments in the purely even case.

The main result of this paper is an explicit construction of
generators of the Feigin--Frenkel center for each simple
Lie algebra $\g=\g_N$ of types $B$, $C$ and $D$.
The construction is based on the properties
of the symmetrizer in the Brauer
algebra which centralizes the action
of the orthogonal or symplectic group
on a tensor product of the vector
representations in the context of the Schur--Weyl duality.
The formulas for the central elements are also presented
in terms of the Howe duality involving a representation of
the Lie algebra $\g_N$ on the symmetric and exterior
powers of the vector space $\CC^N$ and a dual
action of the Lie algebra $\sll_2$.

As a corollary of the main result we get explicit
generators of the respective commutative subalgebras of the
universal enveloping algebras $\U(\g_N[t])$ and the
``shift of argument subalgebras" of $\U(\g_N)$
(see \cite{fft:gm} and \cite{r:si}), which leads
to an explicit construction of higher order Hamiltonians
of the Gaudin model associated with $\g_N$. Moreover,
we use the symmetrizer in the Brauer
algebra to produce commutative
Bethe-type subalgebras in the Drinfeld Yangians $\Y(\g_N)$
and show that their graded images  coincide with the respective
commutative subalgebras of $\U(\g_N[t])$.

To describe the results in more detail, recall that
the affine Kac--Moody algebra $\wh\g$
is the central
extension
\beql{km}
\wh\g=\g\tss[t,t^{-1}]\oplus\CC K,\qquad
\g\tss[t,t^{-1}]=\g\ot\CC[t,t^{-1}],
\eeq
where $\CC[t,t^{-1}]$ is the algebra of Laurent
polynomials in $t$. The vacuum module $V_{-h^{\vee}}(\g)$
at the critical level over $\wh\g$ is defined as
the quotient of the
universal enveloping algebra $\U(\wh\g)$ by the left
ideal generated by $\g[t]$ and $K+h^{\vee}$.
The center of the vertex algebra $V_{-h^{\vee}}(\g)$
is defined by
\ben
\z(\wh\g)
=\{S\in V_{-h^{\vee}}(\g)\ |\ \g[t]\ts S=0\}.
\een
Any element of $\z(\wh\g)$
is called a {\it Segal--Sugawara vector\/}.
The vertex algebra axioms imply that the center is a commutative
associative algebra.
As a vector space, the vacuum module $V_{-h^{\vee}}(\g)$
can be identified with the universal enveloping algebra
$\U(\g_-)$ of the Lie algebra
$\g_-=t^{-1}\g[t^{-1}]$. Moreover, there is
an injective homomorphism $\z(\wh\g)\hra\U(\g_-)$ so that
$\z(\wh\g)$ can be viewed as
a commutative subalgebra of $\U(\g_-)$; see e.g. \cite[Sec.~3.3]{f:lc}.
For any element $S$ of $\U(\g_-)$ denote by
$\overline S$ its {\it symbol\/}, i.e.,
the image in the associated graded algebra
$\text{gr}\ts\U(\g_-)\cong \Sr(\g_-)$.
A set of elements
$
S_1,\dots,S_n\in \z(\wh\g),
$
$n=\text{rank}\ts\g$,
is called a {\it complete set of Segal--Sugawara vectors\/} if
the corresponding symbols
$\overline S_1,\dots,\overline S_n$ coincide with
the images of certain algebraically independent
generators of the algebra
of invariants $\Sr(\g)^{\g}$ under the embedding
$\Sr(\g)\hookrightarrow \Sr(\g_-)$ defined by
the assignment $x\mapsto x\ts t^{-1}$ for $x\in\g$.
According to Feigin and Frenkel~\cite{ff:ak}, a complete set of
Segal--Sugawara vectors exists for any $\g$, and
$\z(\wh\g)$ is the algebra of polynomials in infinitely
many variables,
\beql{genz}
\z(\wh\g)=\CC[T^{\tss r}S_l\ |\ l=1,\dots,n,\ \ r\geqslant 0],
\eeq
where $T$ denotes the translation operator on the vertex
algebra; see Sec.~\ref{sec:pva}.

Now let
$\g=\g_N$ be the orthogonal Lie algebra $\oa_N$ or
the symplectic Lie algebra $\spa_N$ associated
with a nondegenerate bilinear form on $\CC^N$
defined by the respective symmetric
or skew-symmetric matrix $G$.
Recall that the degrees of algebraically independent
generators of the algebra $\Sr(\g_N)^{\g_N}$
are $2,4,\dots,2n$ for $\g_N=\oa_{2n+1}$ and $\g_N=\spa_{2n}$
(types $B$ and $C$, respectively), while for $\g_N=\oa_{2n}$ (type $D$)
the degrees are $2,4,\dots,2n-2,n$.
Introduce standard generators $F_{ij}$ of $\g_N$;
see Sec.~\ref{sec:pva} for the definitions.
Then
$F_{ij}[r]=F_{ij}\ts t^r$ with $r\in\ZZ$ are the corresponding
elements of the loop algebra $\g_N[t,t^{-1}]$. We will need
the extended Lie algebra $\wh\g_N\oplus\CC\tau$ with
the element $\tau$ satisfying the commutation relations
\ben
\big[\tau,F_{ij}[r]\tss\big]=-r\ts F_{ij}[r-1],\qquad
\big[\tau,K\big]=0.
\een
Given $r\in\ZZ$, combine generators of the loop
algebra into the matrix $F[r]=\big[F_{ij}[r]\big]$. Furthermore,
introduce the matrix $\Phi=[\Phi_{ij}]$ by $\Phi=\tau+F[-1]$ so that
$\Phi_{ij}=\de_{ij}\tau+F_{ij}[-1]$,
and identify it with the element
\ben
\Phi=\sum_{i,j=1}^N e_{ij}\ot \Phi_{ij}\in\End\CC^{N}\ot\U,
\een
where the $e_{ij}$ are the standard matrix units and $\U$ stands
for the universal enveloping algebra of
$\wh\g_N\oplus\CC\tau$.
For each $a\in\{1,\dots,m\}$
the element $\Phi_a$ of the algebra
\beql{tenprk}
\underbrace{\End\CC^{N}\ot\dots\ot\End\CC^{N}}_m{}\ot\U
\eeq
is defined by
\ben
\Phi_a=\sum_{i,j=1}^{N}
1^{\ot(a-1)}\ot e_{ij}\ot 1^{\ot(m-a)}\ot \Phi_{ij}.
\een

Consider the {\it Brauer algebra\/} $\Bc_m(\om)$
over $\CC(\om)$ and let $S^{(m)}\in\Bc_m(\om)$
denote the {\it symmetrizer\/} in $\Bc_m(\om)$ which is
the primitive idempotent associated with the trivial representation
of $\Bc_m(\om)$. A few equivalent expressions
for the symmetrizer are given in Sec.~\ref{subsec:sba}.
The classical orthogonal and symplectic groups
$\Or_N$ and $\Spr_N$
act naturally on the tensor product space $(\CC^N)^{\ot m}$.
By the results of
Brauer~\cite{b:aw}, the centralizer of this action
in the endomorphism algebra of the tensor product space
coincides with the homomorphic
image of the algebra $\Bc_m(\om)$
with the parameter $\om$ specialized to $N$ and $-N$, respectively,
in the orthogonal and symplectic case.
Using the natural action of the Brauer
algebra on $(\CC^{N})^{\ot m}$
we will regard the symmetrizer $S^{(m)}$
as an element of the algebra \eqref{tenprk} with
the identity component in $\U$. In the orthogonal case,
taking trace over
all $m$ copies of $\End\CC^{N}$, write the following
element as a polynomial
in $\tau$ with coefficients in the universal enveloping
algebra $\U(t^{-1}\g_N[t^{-1}])$:
\beql{deftr}
\tr\ts S^{(m)} \Phi_1\dots \Phi_m
=\phi^{}_{m\tss0}\ts\tau^m+\phi^{}_{m\tss1}\ts\tau^{m-1}
+\dots+\phi^{}_{m\tss m}.
\eeq
If $\g_N=\oa_{2n}$, we also define the element
$\phi^{\tss\prime}_n=\Pf \wt F[-1]$ as
the (noncommutative) Pfaffian of the skew-symmetric matrix $\wt F[-1]=F[-1]\tss G$,
\beql{genpf}
\Pf \wt F[-1]=\frac{1}{2^nn!}\sum_{\si\in\Sym_{2n}}\sgn\si\cdot
\wt F_{\si(1)\ts\si(2)}[-1]\dots \wt F_{\si(2n-1)\ts\si(2n)}[-1].
\eeq

Defining analogues of the elements $\phi^{}_{m\tss k}$
in the symplectic case requires an extra care.
The image of the element
$S^{(m)}\in\Bc_m(-2n)$ in \eqref{tenprk} with $N=2n$
is well-defined for $m\leqslant n+1$ but it is zero
for $m=n+1$, and the specialization of
$S^{(m)}\in\Bc_m(\om)$ at $\om=-2n$
is not defined for
$n+2\leqslant m\leqslant 2n$. However, using a symmetry
property of the matrix $\Phi$, we may regard the product
$\Phi_1\dots \Phi_m$ as a polynomial in $\tau$ with
coefficients in the algebra of the form \eqref{tenprk} where
each factor $\End\CC^{2n}$ is replaced by the direct sum of the subspaces
of symplectic matrices and scalar matrices.
Restricting the operator $S^{(m)}$
to the tensor products of these subspaces allows us to
write
\beql{deftrsp}
\frac{1}{n-m+1}\ts
\tr\ts S^{(m)} \Phi_1\dots \Phi_m
=\phi^{}_{m\tss0}\ts\tau^m+\phi^{}_{m\tss1}\ts\tau^{m-1}
+\dots+\phi^{}_{m\tss m},
\eeq
so that the elements $\phi^{}_{m\tss k}$ are well-defined
for all values of the parameters with $1\leqslant m\leqslant 2n$.

\begin{mthm} {\rm(i)}\quad All elements $\phi^{}_{m\tss k}$ are
Segal--Sugawara vectors for $\g_N$.\par
\noindent
{\rm(ii)}\quad
$\{\phi^{}_{2\tss 2},\phi^{}_{4\tss 4},\dots,
\phi^{}_{2n\ts 2n}\}$
is a complete set of Segal--Sugawara vectors
for $\oa_{2n+1}$ and $\spa_{2n}$.\par
\noindent
{\rm(iii)}\quad
$\{\phi^{}_{2\tss 2},\phi^{}_{4\tss 4},\dots,\phi^{}_{2n-2\ts 2n-2},
\phi^{\tss\prime}_n\}$
is a complete set of Segal--Sugawara vectors
for $\oa_{2n}$.
\end{mthm}

The Segal--Sugawara vectors $\phi^{}_{m\tss k}$ defined by
\eqref{deftr} and \eqref{deftrsp}
admit an alternative presentation with the use
of the respective Howe dual pair $(\sll_2,\g_N)$; see \cite{h:pi}.
Namely, the operator $S^{(m)}$ projects
the vector space $(\CC^N)^{\ot m}$ to a subspace
of symmetric tensors in the orthogonal case and to
a subspace of skew-symmetric tensors in the symplectic case.
Using natural identifications of these subspaces
with the respective homogeneous components of the
symmetric algebra $\Sr(\CC^N)$ and exterior algebra
$\La(\CC^{N})$, we find that the action of $S^{(m)}$
coincides with that of the extremal projector $p=p(\sll_2)$
associated with the commuting action of the Lie
algebra $\sll_2$. Therefore, the traces on the left hand sides
of \eqref{deftr} and \eqref{deftrsp} can be replaced by the
respective traces $\tr\ts p\ts\Phi^{(m)}$ taken over
the subspace of $\sll_2$-singular vectors
in the homogeneous components of degree $m$ in
$\Sr(\CC^N)$ and $\La(\CC^{N})$, where $\Phi^{(m)}$
denotes the restriction of the element $\Phi_1\dots\Phi_m$
to this subspace.

Note that together with the formulas of \cite{cm:ho}, \cite{ct:qs}
for the Segal--Sugawara vectors in type $A$, the Main Theorem leads to
a new proof of the Feigin--Frenkel theorem for all classical types;
see \cite[Ch.~3]{f:lc} for a detailed discussion of this explicit approach
exemplified by the case of $\sll_2$.

\medskip

I am grateful to Alexander Chervov and Eric Ragoucy
for collaboration in \cite{cm:ho} and \cite{mr:mm}
and discussions of examples and approaches to the Sugawara operators
for the classical Lie algebras. I also thank Tomoyuki Arakawa,
Joe Chuang, Pavel Etingof and Evgeny Mukhin
for stimulating discussions.

\section{Affine vertex algebras in types $B$, $C$ and $D$}
\label{sec:pva}
\setcounter{equation}{0}

We will be working with the affine vertex
algebras $V_{\ka}(\g_N)$ associated
with the classical Lie algebras $\g_N$. To introduce
necessary definitions, consider the vector space $\CC^N$
with its canonical basis $e_1,\dots,e_N$ and
equip it with a nondegenerate symmetric
or skew-symmetric bilinear form
\beql{formg}
\langle e_i,e_j\rangle=g_{ij}
\eeq
so that the matrix $G=[g_{ij}]$ is
symmetric or skew-symmetric, respectively,
$G^{\tss t}={\pm}\ts G$, where $t$ denotes the standard
matrix transposition.
The skew-symmetric case may only occur
for even $N$. The respective classical
group $\Gr_N=\Or_N$ or $\Gr_N=\Spr_{N}$ is defined as the group
of complex matrices
preserving the form \eqref{formg}.
For any $N\times N$ matrix $A$ set
$
A'=G A^t\ts G^{-1}.
$
For matrices with entries in $\CC$ this is just the
transpose with respect to the form determined
by the inverse matrix $G^{-1}=[\overline{g}_{ij}]$.
Denote by $E_{ij}$, $1\leqslant i,j\leqslant N$,
the standard basis vectors of the Lie
algebra $\gl_N$.
Introduce the elements $F_{ij}$
of $\gl_N$ by the formulas
\beql{fijelem}
F_{ij}=E_{ij}-\sum_{k,\tss l=1}^N
g^{}_{ik}\tss \overline g^{}_{lj}\tss E_{lk}.
\eeq
The Lie subalgebra
of $\gl_N$ spanned by the elements $F_{ij}$ is isomorphic to
the orthogonal Lie algebra $\oa_N$ associated with $\Or_N$
in the symmetric case
and to the symplectic Lie algebra $\spa_N$
associated with $\Spr_N$
in the skew-symmetric case. Denote by $F$ the $N\times N$ matrix
whose $(i,j)$ entry is $F_{ij}$. We shall also regard $F$ as the
element
\ben
F=\sum_{i,j=1}^N e_{ij}\ot F_{ij}\in \End\CC^N\ot \U(\g_N).
\een
The definition \eqref{fijelem} can be written
in the matrix form as $F=E-E^{\tss\prime}$,
where $E=[E_{ij}]$.
Consider the permutation operator
\beql{pmatrix}
P=\sum_{i,j=1}^N e_{ij}\ot e_{ji}\in \End \CC^N\ot\End \CC^N
\eeq
and set
\beql{rprimeq}
Q=G_1\tss P^{\tss t}\ts G_1^{-1}=G_2\tss P^{\tss t}\ts G_2^{-1},
\qquad P^{\tss t}=\sum_{i,j=1}^N e_{ij}\ot e_{ij},
\eeq
where $G_1=G\ot 1$ and $G_2=1\ot G$.
Note that the operators $P$ and $Q$ satisfy
the relations
\beql{pq}
P^2=1,\qquad Q^2=N\ts Q,\qquad PQ=QP=\begin{cases}
\phantom{-}Q\quad&\text{in the symmetric case,}\\
-Q\quad&\text{in the skew-symmetric case.}
\end{cases}
\eeq
The defining relations of the algebra $\U(\g_N)$
have the form
\beql{deff}
F_1\ts F_2-F_2\ts F_1=(P-Q)\ts F_2-F_2\ts (P-Q)
\eeq
together with the relation $F+F^{\tss\prime}=0$,
where both sides in \eqref{deff}
are regarded as elements of the algebra
$\End\CC^N\ot\End\CC^N\ot \U(\g_N)$ and
\beql{eonetwo}
F_1=\sum_{i,j=1}^N e_{ij}\ot 1\ot F_{ij},\qquad
F_2=\sum_{i,j=1}^N 1\ot e_{ij}\ot F_{ij}.
\eeq

Now consider the affine
Kac--Moody algebra $\wh\g_N=\g_N\tss[t,t^{-1}]\oplus\CC K$
and set $F_{ij}[r]=F_{ij}t^r$ for any $r\in\ZZ$.
As with the Lie algebra $\g_N$, regard the matrix $F[r]=\big[F_{ij}[r]\big]$
as the element
\ben
F[r]=\sum_{i,j=1}^N e_{ij}\ot F_{ij}[r]\in \End\CC^N\ot \U(\wh\g_N).
\een
By analogy with \eqref{deff},
the defining relations of the algebra $\U(\wh\g_N)$
can be written in the form (cf. \cite{k:id}):
\begin{multline}\label{deffaff}
F[r]_1\ts F[s]_2-F[s]_2\ts F[r]_1=(P-Q)\ts F[r+s]_2-F[r+s]_2\ts (P-Q)\\
{}+{}\begin{cases}
r\tss\de_{r,-s}\tss(P-Q)\tss K\quad&\text{in the orthogonal case,}\\
2\tss r\tss\de_{r,-s}\tss(P-Q)
\tss K\quad&\text{in the symplectic case.}
\end{cases}
\end{multline}

A {\it vertex algebra\/} $V$ is a vector space with
the additional data $(Y, T, \vac)$, where
the state-field correspondence $Y$ is a map
$Y: V \rightarrow \End V[[z,z^{-1}]]$,
the infinitesimal translation $T$ is an operator $T: V \rightarrow V$,
and $\vac \in V$ is a vacuum vector.
For the axioms satisfied by these data
see e.g., \cite{dk:fa},
\cite{fb:va}, \cite{k:va}.
For any element $a \in V$ we write
\ben
Y(a,z) = \sum_{n\in\ZZ} a_{(n)} z^{-n-1}, \quad
a_{(n)} \in \End V.
\een
The {\it center\/} of a vertex algebra $V$ is its commutative
vertex subalgebra spanned by all vectors $b\in V$
such that $a_{(n)}\ts b=0$ for all $a\in V$ and $n\geqslant 0$.

For any $\ka\in\CC$
define the {\it affine vertex algebra\/}
$V_{\ka}(\g_N)$ as the quotient of the
universal enveloping algebra $\U(\wh\g_N)$ by the left
ideal generated by $\g_N[t]$ and $K-\ka$.
The vacuum vector is $1$ and
the translation operator is determined by
\ben
T:1\mapsto 0, \quad\big[T, K\big]=0\fand
\big[T, F_{ij}[r]\tss\big]= -r\tss F_{ij}[r-1].
\een
The state-field correspondence $Y$
is defined by setting $Y(1,z)=\text{id}$, $Y(K,z)=K$,
\beql{basfi}
Y(F_{ij}[-1],z)=F_{ij}(z):=
\sum_{r\in\ZZ} F_{ij}[r]\tss z^{-r-1},
\eeq
and then extending the map to the whole of $V_{\ka}(\g_N)$
with the use of normal ordering; see
{\it loc. cit.}

\section{Brauer algebra and Howe pairs}
\label{sec:bra}
\setcounter{equation}{0}

Let $\om$ be an indeterminate and
$m$ a positive integer.
The {\it Brauer algebra\/} $\Bc_m(\om)$ over the field
$\CC(\om)$ is generated by
$2m-2$ elements
$s_1,\dots,s_{m-1},\ep_1,\dots,\ep_{m-1}$
subject to the defining relations
\beql{bradr}
\begin{aligned}
s_i^2&=1,\qquad \ep_i^2=\om\ts \ep_i,\qquad
\su_i\ee_i=\ee_i\su_i=\ee_i,\qquad i=1,\dots,m-1,\\
s_i s_j &=s_js_i,\qquad \ee_i \ee_j = \ee_j \ee_i,\qquad
\su_i \ee_j = \ee_j\su_i,\qquad
|i-j|>1,\\
s_is_{i+1}s_i&=s_{i+1}s_is_{i+1},\qquad
\ee_i\ee_{i+1}\ee_i=\ee_i,\qquad \ee_{i+1}\ee_i\ee_{i+1}=\ee_{i+1},\\
\su_i\ee_{i+1}\ee_i&=\su_{i+1}\ee_i,\qquad \ee_{i+1}\ee_i
\su_{i+1}=\ee_{i+1}\su_i,\qquad
i=1,\dots,m-2.
\end{aligned}
\non
\end{equation}
This algebra was originally defined by
Brauer~\cite{b:aw} in terms of
diagrams and it is well-known it admits
the above presentation.
The subalgebra of $\Bc_m(\om)$ generated over $\CC$
by $s_1,\dots,s_{m-1}$
is isomorphic to the group algebra $\CC[\Sym_m]$ so that $s_i$
can be identified with the transposition $(i\ i+1)$.
For any $1\leqslant i<j\leqslant n$
the transposition $s_{ij}=(i\ts j)$ can be regarded
as an element of $\Bc_m(\om)$. We denote by $\ep_{ij}$
the elements of $\Bc_m(\om)$ defined by
$\ep_{j-1\ts j}=\ep_{j-1}$ and
$\ep_{i\ts j}=s_{i\ts j-1}\tss \ep_{j-1}\tss s_{i\ts j-1}$ for $i<j-1$.
The Brauer algebra $\Bc_{m-1}(\om)$ can be regarded
as a natural subalgebra of $\Bc_m(\om)$.

\subsection{Symmetrizer in the Brauer algebra}
\label{subsec:sba}

Consider the one-dimensional representation of $\Bc_m(\om)$
where each element $s_{ij}$ acts as the identity operator and each
element $\ep_{ij}$ acts as the zero operator. The corresponding
idempotent $S^{(m)}\in\Bc_m(\om)$ then has the properties
\beql{mulse}
s_{ij}\tss S^{(m)}=S^{(m)}\tss s_{ij}=S^{(m)}\Fand
\ep_{ij}\tss S^{(m)}=S^{(m)}\tss\ep_{ij}=0
\eeq
for all $i,j$, and it
can be given by several
equivalent formulas. The first one relies on the
Jucys--Murphy elements for the Brauer algebra
and their eigenvalues in the irreducible representations
found independently in \cite{lr:rh} and \cite{n:yo}.
The corresponding formula reads
\beql{symjm}
S^{(m)}=\prod_{r=2}^m\frac{1}{r\tss(\om+2r-4)}
\Big(1+\sum_{i=1}^{r-1}(s_{i\tss r}-\ep_{i\tss r})\Big)
\Big(\om+r-3+\sum_{i=1}^{r-1}(s_{i\tss r}-\ep_{i\tss r})\Big).
\eeq
The ordering of the factors is irrelevant as they commute pairwise.
Two more expressions for the symmetrizer are related to
{\it fusion procedures\/} for the Brauer algebra;
see \cite{im:fp}, \cite{imo:fp} and references therein.
We have
\beql{symfp}
S^{(m)}=\frac{1}{m!}\prod_{1\leqslant i<j\leqslant m}
\Big(1-\frac{\ep_{ij}}{\om+i+j-3}\Big)
\prod_{1\leqslant i<j\leqslant m}
\Big(1+\frac{s_{ij}}{j-i}\Big)
\eeq
and
\beql{symfpnew}
S^{(m)}=\frac{1}{m!}
\prod_{1\leqslant i<j\leqslant m}
\Big(1+\frac{s_{ij}}{j-i}-\frac{\ep_{ij}}
{\om/2+j-i-1}\Big),
\eeq
where all products are taken in the lexicographic order
on the pairs $(i,j)$. It is well-known by Jucys~\cite{j:yo} that
the second product in \eqref{symfp} coincides
with the symmetrizer in the group algebra $\CC[\Sym_m]$:
\ben
H^{(m)}:=\frac{1}{m!}\prod_{1\leqslant i<j\leqslant m}
\Big(1+\frac{s_{ij}}{j-i}\Big)=\frac{1}{m!}\sum_{s\in \Sym_m} s.
\een
We will need yet another formula for $S^{(m)}$.

\bpr\label{prop:symexp}
We have
\beql{symexp}
S^{(m)}=H^{(m)}\ts\sum_{r=0}^{\lfloor m/2\rfloor}
\frac{(-1)^r}{2^{\tss r}\ts r!}\binom{\om/2+m-2}{r}^{-1}
\sum_{i_a<j_a}\ep_{i_1j_1}\ep_{i_2j_2}\dots \ep_{i_rj_r}
\eeq
with the second sum taken over the {\rm(}unordered\tss{\rm)}
sets of disjoint pairs
$\{(i_1,j_1),\dots,(i_r,j_r)\}$ of indices from $\{1,\dots,m\}$.
\epr

\bpf
Note that for each $r$ the second sum commutes with
any element of $\Sym_m$ and hence commutes with $H^{(m)}$.
Now we apply \eqref{symfp}
and write $H^{(m)}=H^{(m-1)}\tss H^{(m)}$.
Since
$
\ep_{im}\tss\ep_{jm}\tss H^{(m-1)} =\ep_{im}\tss s_{ij}\tss H^{(m-1)}
=\ep_{im}\tss H^{(m-1)}
$
for distinct $i$ and $j$, we get
\ben
\Big(1-\frac{\ep_{1m}}{\om+m-2}\Big)\dots
\Big(1-\frac{\ep_{m-1\ts m}}{\om+2m-4}\Big)\tss H^{(m-1)}
=\Big(1-\frac{\ep_{1m}+\dots+\ep_{m-1\ts m}}{\om+2m-4}\Big)\tss H^{(m-1)}
\een
via a straightforward calculation.
Using induction on $m$ we come to verifying that
if \eqref{symexp} holds for $m$ replaced by $m-1$, then
\ben
S^{(m-1)}\ts \Big(1-\frac{\ep_{1m}+\dots+\ep_{m-1\ts m}}
{\om+2m-4}\Big)\tss H^{(m)}
\een
will coincide with the right hand side of \eqref{symexp}.
If none of the pairs $(i_a,j_a)$ contains
an index $l\in\{1,\dots,m-1\}$, then
the product $\ep_{i_1j_1}\ep_{i_2j_2}\dots \ep_{i_rj_r}\ts\ep_{l\tss m}$
of $r+1$ factors
will contribute to the respective sum on
the right hand side of \eqref{symexp}.
Otherwise,
if $l=i_a$ or $l=j_a$ for some $a$, then
$\ep_{i_1j_1}\ep_{i_2j_2}\dots \ep_{i_rj_r}\ts\ep_{l\tss m}\tss H^{(m)}
=\ep_{i_1j_1}\ep_{i_2j_2}\dots \ep_{i_rj_r}\tss H^{(m)}$,
and the result follows by combining similar terms.
\epf

\subsection{Dual pairs and extremal projector}
\label{subsec:sv}

Consider the natural right action of the classical group
$\Gr_N=\Or_N$ or $\Gr_N=\Spr_N$
in the
tensor product space
\beql{cntenpr}
\underbrace{\CC^N\ot\CC^N\ot\dots\ot\CC^N}_m.
\eeq
Label the copies of the vector space $\CC^N$ with the
numbers $1,\dots,m$ and for each pair
of indices $1\leqslant i<j\leqslant m$
denote by $P_{ij}$ and $Q_{ij}$ the operators on
the vector space \eqref{cntenpr} which act as
the respective operators $P$ and $Q$ defined in \eqref{pmatrix}
and \eqref{rprimeq} on the tensor product on the $i$-th and $j$-th
copies of $\CC^N$ and act as the identity operators on each
of the remaining copies. The commuting action
of the Brauer algebra $\Bc_m(N)$ in the orthogonal case
is given by the assignment
\beql{braact}
s_{ij}\mapsto P_{ij},\qquad \ep_{ij}\mapsto Q_{ij},\qquad
i< j;
\eeq
see \cite{b:aw}. Similarly, in the symplectic case with $N=2n$
the commuting
action of the Brauer algebra $\Bc_m(-2n)$ on the space
\eqref{cntenpr}
is given by
\beql{braactsy}
s_{ij}\mapsto -\ts P_{ij},\qquad \ep_{ij}\mapsto -\ts Q_{ij},\qquad
i< j.
\eeq

We will need particular
symmetric and skew-symmetric forms
\eqref{formg}. We will be using the involution on the set
of indices $\{1,\dots,N\}$ defined by $i\mapsto i^{\tss\prime}=N-i+1$
and take $G=[g_{ij}]$ with
\beql{formgss}
g_{ij}=\begin{cases}
\phantom{\ve_i\ts}\de_{i\ts j^{\tss\prime}}\quad&
\text{in the symmetric case,}\\
\ve_i\ts\de_{i\ts j^{\tss\prime}}\quad&\text{in the skew-symmetric case,}
\end{cases}
\eeq
where $\ve_i=1$ for $i=1,\dots,n$ and
$\ve_i=-1$ for $i=n+1,\dots,2n$.

In the orthogonal case the element $H^{(m)}\in\Bc_m(N)$ acts as
the symmetrizer on the space \eqref{cntenpr}.
Denote by $\Pc_N$
the space of polynomials
in variables $z_1,\dots,z_N$ and
identify the image $H^{(m)}(\CC^N)^{\ot m}$
with the subspace $\Pc_N^{\tss m}$
of homogeneous polynomials of degree $m$
via the isomorphism
\beql{symsub}
H^{(m)}(e_{i_1}\ot\dots\ot e_{i_m})\mapsto z_{i_1}\dots z_{i_m}.
\eeq
Proposition~\ref{prop:symexp} allows us
to consider $S^{(m)}$ as an operator
on the space $\Pc_N^{\tss m}$.
This operator commutes
with the action of the orthogonal group $\Or_N$ on this space and so
it can be expressed in terms of the
action of the generators of the Lie algebra $\sll_2$
provided by the Howe duality~\cite[Sec.~3.4]{h:pi}.
Let $\{e,f,h\}$ be the standard basis of $\sll_2$.
The action of the basis elements
on the space $\Pc_N$
is given by
\ben
e\mapsto -\frac12\ts\sum_{i=1}^N \di_{i}\tss
\di_{\tss i^{\tss\prime}},\qquad
f\mapsto \frac12\ts\sum_{i=1}^N z_i\tss z_{\tss {i\tss}'},\qquad
h\mapsto -\frac{N}{2}-\sum_{i=1}^N
z_i\ts\di_i,
\een
where $\di_i$ denotes the partial derivative over $z_i$.
Recall that the {\it extremal projector\/} $p$
for the Lie algebra $\sll_2$ (see \cite{ast:po})
is given by the formula
\beql{ep}
p=1+\sum_{r=1}^{\infty}
\ts\frac{(-1)^r}{r!\ts(h+2)\dots
(h+r+1)}\tss f^{\tss r} e^r.
\eeq
The projector $p$ is an element
of an extension of the universal enveloping algebra $\U(\sll_2)$
and it possesses the properties
$e\tss p=p\tss f=0$. The basis
element $h$ of $\sll_2$ acts on the space $\Pc_N^{\tss m}$
as multiplication by the scalar $-N/2-m$. Furthermore,
this space is annihilated by the action of $e^r$ with $r>m/2$.
Hence, using \eqref{ep} we may regard
$p$ as an operator on $\Pc_N$.
The image of this operator coincides with the
subspace of $\sll_2$-singular vectors.

On the other hand, for each value of
$r=0,\dots,\lfloor m/2\rfloor$
the image of the second sum in
\eqref{symexp} under the action \eqref{braact}
coincides with
the operator
\beql{qop}
\sum Q_{i_1j_1}Q_{i_2j_2}\dots Q_{i_rj_r}
\eeq
which preserves
the subspace $H^{(m)}(\CC^N)^{\ot m}$. Regarding
\eqref{qop} as an operator on $\Pc_N^{\tss m}$
via the isomorphism \eqref{symsub}, we verify that
it coincides with the action of the element
$(-2)^{\tss r} f^r e^r/r!$. Comparing \eqref{symexp} and \eqref{ep}
we come to the following result.

\bpr\label{prop:expsm}
The operator $S^{(m)}$ on the space $\Pc_N^{\tss m}$
coincides with the restriction of the action of
the extremal projector $p$ to this subspace.
\qed
\epr

In the symplectic case
the element $H^{(m)}\in\Bc_m(-2n)$ acts as
the anti-symmetrizer on \eqref{cntenpr}.
We will regard the exterior algebra $\La_{2n}=\La(\CC^{2n})$
as the space of polynomials in the anti-commuting
variables $\ze_1,\dots,\ze_{2n}$ and denote by
$\La^m_{2n}$ the subspace of homogeneous polynomials
of degree $m$. We
identify the image $H^{(m)}(\CC^{2n})^{\ot m}$
with $\La^m_{2n}$ via the isomorphism
\beql{symsubex}
H^{(m)}(e_{i_1}\ot\dots\ot e_{i_m})\mapsto
\ze_{i_1}\wedge\dots\wedge \ze_{i_m}.
\eeq
If the condition $m\leqslant n$ holds, then $S^{(m)}$
is a well-defined operator on $\La^m_{2n}$
by Proposition~\ref{prop:symexp}.
This operator commutes
with the action of the symplectic group $\Spr_{2n}$ on this space
and can be expressed in terms of the
action of the generators of the Lie algebra $\sll_2$
provided by the skew Howe duality~\cite[Sec.~3.8]{h:pi}.
The action of the basis elements of $\sll_2$
on the space $\La_{2n}$
is given by
\ben
e\mapsto \sum_{i=1}^n \di_{i}\wedge \di_{\tss i^{\tss\prime}},\qquad
f\mapsto -\sum_{i=1}^n \ze_i\wedge \ze_{\tss i^{\tss\prime}},\qquad
h\mapsto n-\sum_{i=1}^n
\big(\ze_i\wedge \di_i+\ze_{i^{\tss\prime}}
\wedge \di_{\tss i^{\tss\prime}}\big),
\een
where $\di_i$ denotes the (left) partial derivative over $\ze_i$
and $i^{\tss\prime}=2n-i+1$. Note that the basis element $h$
acts as multiplication by the scalar $n-m$ on $\La_{2n}^m$.
Since the representation $\La_{2n}$ of $\sll_2$ is finite-dimensional,
the subspace $\La_{2n}^m$ can contain singular vectors only if
$n-m\geqslant 0$ which is consistent with our condition
on the parameters. Under this condition
the extremal projector $p$ defined in \eqref{ep}
can be regarded as an operator on $\La_{2n}^m$
projecting to a subspace of singular vectors.
Using again Proposition~\ref{prop:symexp} we come
to the following.

\bpr\label{prop:expsmsp}
If $m\leqslant n$ then
the operator $S^{(m)}$ on the space $\La_{2n}^{m}$
coincides with the restriction of the action of
the extremal projector $p$ to this subspace.
\qed
\epr

\subsection{Noncommutative characteristic maps}
\label{subsec:cm}

Extend the notation
\eqref{eonetwo} to multiple tensor products
\beql{tenprku}
\underbrace{\End\CC^{N}\ot\dots\ot\End\CC^{N}}_m{}\ot\U(\g_N),
\eeq
and for an arbitrary element $C$ of the respective
Brauer algebra $\Bc_m(N)$ or $\Bc_m(-N)$ consider
the element $\tr\ts C\tss F_1\dots F_m$ of $\U(\g_N)$,
where the trace is taken over all $m$ copies of $\End\CC^{N}$,
and the element $C$ is identified with its image in the
algebra \eqref{tenprku} acting as the identity operator
on $\U(\g_N)$. Observe that $\tr\ts C\tss F_1\dots F_m$
belongs to the subalgebra $\U(\g_N)^{\Gr_N}$
of $\Gr_N$-invariants in the universal enveloping algebra; this
follows by noting that for any matrix $\hb$ such that $\hb^t\in \Gr_N$ we have
\begin{multline}\non
\tr\ts C\ts \hb_1 F_1\hb_1^{-1}\dots \hb_m F_m \hb_m^{-1}
=\tr\ts C\ts \hb_1\dots \hb_m
\tss F_1\dots F_m \hb_1^{-1}\dots\hb_m^{-1}\\
{}=\tr\ts \hb_1^{-1}\dots\hb_m^{-1}
\tss C\ts \hb_1\dots \hb_m
\tss F_1\dots F_m=\tr\ts C\tss F_1\dots F_m,
\end{multline}
where we used the cyclic property of trace and the fact
that the action of the element $C$ commutes with the right action
of $\Gr_N$. The product $\hb_1\dots \hb_m$ is the image of $\hb^t$ under this
action of $\Gr_N$.
On the other hand, the algebra of invariants
$\U(\g_N)^{\Gr_N}$ is isomorphic to the algebra of symmetric
polynomials $\CC[l_1^2,\dots,l_n^2]^{\Sym_n}$ in $n$ variables,
if $N=2n+1$ or $N=2n$ via the Harish-Chandra isomorphism;
see e.g. \cite[Sec.~7.4]{d:ae}. Here we use
basis elements $h_1,\dots,h_n$ of a Cartan subalgebra
$\h_N$ of $\g_N$ and $l_i=h_i+\rho_i$
for $i=1,\dots,n$ are their shifts by constants $\rho_i$
determined by the half-sum of the positive roots.
This gives (noncommutative) Brauer algebra versions
of the well-known characteristic map associated with the
symmetric groups:
\beql{charmap}
\ch:\Bc_m(\om)\to \CC[l_1^2,\dots,l_n^2]^{\Sym_n},\qquad \om=N,-N.
\eeq
We will need to calculate the leading terms of the symmetric
polynomials $\ch(S^{(m)})$. To this end, we can work
with the corresponding commutative version of
the characteristic map \eqref{charmap}
by replacing the algebra of invariants $\U(\g_N)^{\Gr_N}$
with the associated graded algebra $\Sr(\g_N)^{\Gr_N}$ which is
isomorphic to the algebra of the symmetric
polynomials $\CC[h_1^2,\dots,h_n^2]^{\Sym_n}$; see
\cite[Sec.~7.3]{d:ae}. Furthermore, we fix
the respective bilinear form on $\CC^N$ defined in \eqref{formgss}.
The leading term
of the symmetric polynomial $\ch(S^{(m)})$
(with the condition $m\leqslant n$ in the
symplectic case) can now be found as
the trace
\beql{leadterm}
\tr\ts S^{(m)} X_1\dots X_m,
\eeq
taken over all copies of $\End\CC^N$ in the algebra
\ben
\underbrace{\End\CC^{N}\ot\dots\ot\End\CC^{N}}_m{}\ot\Sr(\h_N),
\een
where $X$ is the diagonal
matrix
\ben
X={\rm diag}[h_1,\dots,h_n,0,-h_n,\dots,-h_1]
\For
X={\rm diag}[h_1,\dots,h_n,-h_n,\dots,-h_1]
\een
for $N=2n+1$ or $N=2n$, respectively, with entries in $\Sr(\h_N)$.

\bpr\label{prop:leadsym}
In the orthogonal case, the trace \eqref{leadterm} is zero
if $m$ is odd, while for even $m=2k$ we have
\ben
\tr\ts S^{(2k)} X_1\dots X_{2k}=\frac{N+4k-2}{N+2k-2}
\sum_{1\leqslant i_1\leqslant\dots\leqslant i_k\leqslant n}
\ts h_{i_1}^2\dots h_{i_k}^2.
\een
\epr

\bpf
We will write the diagonal entries of the matrix $X$
as $h_1,\dots,h_N$ by setting $h_{i^{\tss\prime}}=-h_i$
so that $h_{n+1}=0$ if $N=2n+1$.
The left hand side of \eqref{leadterm} then takes the form
\beql{sumis}
\sum_{i_1,\dots,i_m=1}^N h_{i_1}\dots h_{i_m}\ts a_{i_1,\dots,i_m},
\eeq
where the coefficients $a_{i_1,\dots,i_m}$ are defined
as the diagonal entries in the expansion
\ben
S^{(m)}(e_{i_1}\ot\dots\ot e_{i_m})=
a_{i_1,\dots,i_m}(e_{i_1}\ot\dots\ot e_{i_m})+\dots
\een
in terms of the basis vectors of $(\CC^N)^{\ot m}$.
Note that the values $i_a=n+1$ in the case $N=2n+1$
can be excluded in the sum in \eqref{sumis} as $h_{n+1}=0$.
Now we use Proposition~\ref{prop:expsm}.
We have
\ben
\frac{f^{\tss r}}{r!}\mapsto\sum\frac{(z_1z_{1'})^{a_1}}{a_1!}
\dots\frac{(z_nz_{n'})^{a_n}}{a_n!},\qquad
\frac{(-e)^{\tss r}}{r!}\mapsto\sum\frac{(\di_1\di_{1'})^{a_1}}{a_1!}
\dots\frac{(\di_n\di_{n'})^{a_n}}{a_n!},
\een
where each sum is taken over all tuples of nonnegative
integers $a_i$ such that $a_1+\dots+a_n=r$ and for $N=2n+1$
we restrict the action to the subspace of polynomials of degree $0$
in the variable $z_{n+1}$. For
any monomial
$z_{1\phantom{^\prime}}^{b^{}_{1}}
z_{1'}^{b_{1'}}\dots z_{n\phantom{^\prime}}^{b^{}_n}z_{n'}^{b_{n'}}$
of degree $m$ we then have the diagonal coefficient in the expansion
\ben
\frac{f^{\tss r}}{r!}
\frac{(-e)^{\tss r}}{r!}\ts z_{1\phantom{^\prime}}^{b^{}_{1}}
z_{1'}^{b_{1'}}\dots z_{n\phantom{^\prime}}^{b^{}_n}z_{n'}^{b_{n'}}
=\sum \binom{b_1}{a_1}\binom{b_{1'}}{a_1}\dots
\binom{b_n}{a_n}\binom{b_{n'}}{a_n}
\ts z_{1\phantom{^\prime}}^{b^{}_{1}}
z_{1'}^{b_{1'}}\dots z_{n\phantom{^\prime}}^{b^{}_n}z_{n'}^{b_{n'}}
+\dots
\een
with the sum over all tuples of nonnegative
integers $a_i$ such that $a_1+\dots+a_n=r$.
Since $h_{i^{\tss\prime}}=-h_i$, calculating
the coefficient
of the monomial $h_1^{c_1}\dots h_n^{c_n}$ in \eqref{sumis}
involves evaluating the sums
\ben
\sum_{b_i+b_{i'}=c_i}
(-1)^{b_{i^{\tss\prime}}}\binom{b_i}{a_i}\binom{b_{i'}}{a_i}.
\een
If $c_i$ is odd then the sum is zero, while for $c_i=2\tss d_i$
the sum equals $(-1)^{a_i}{\displaystyle\binom{d_i}{a_i}}$.
This implies the first part of the proposition. Now suppose
that $m=2k$ is even. Then the coefficient of the monomial
$h_1^{2d_1}\dots h_n^{2d_n}$ in \eqref{sumis}
equals
\ben
\bal
\sum_{r=0}^{k}
(-1)^r\binom{N/2+2k-2}{r}^{-1}&
\sum (-1)^{a_1}\binom{d_1}{a_1}\dots (-1)^{a_n}\binom{d_n}{a_n}\\
{}&=\sum_{r=0}^{k}
\binom{N/2+2k-2}{r}^{-1}\binom{k}{r}=\frac{N/2+2k-1}{N/2+k-1}
\eal
\een
thus completing the proof.
\epf

\bpr\label{prop:leadsksym}
In the symplectic case with $m\leqslant n$,
the trace \eqref{leadterm} is zero
if $m$ is odd, while for even values $m=2k$
we have
\ben
\tr\ts S^{(2k)} X_1\dots X_{2k}=(-1)^k\ts\frac{n-2k+1}{n-k+1}
\sum_{1\leqslant i_1<\dots< i_k\leqslant n}
\ts h_{i_1}^2\dots h_{i_k}^2.
\een
\epr

\bpf
Now we use the formulas
\eqref{braactsy} for the action of the Brauer algebra $\Bc_m(-N)$
with $N=2n$ to calculate the corresponding sum \eqref{sumis}.
Since $H^{(m)}$ now acts as the anti-symmetrizer,
the condition
$m\leqslant n$ implies that
a basis in
the vector space $H^{(m)}(\CC^N)^{\ot m}$ is formed
by the vectors
\ben
H^{(m)}(e_{c^{}_1}\ot e_{c^{\tss\prime}_1}\ot\dots\ot
e_{c^{}_p}\ot e_{c^{\tss\prime}_p}
\ot e_{b^{}_1}\ot\dots\ot e_{b^{}_{m-2p}}),
\qquad p=0,1,\dots,\lfloor m/2\rfloor,
\een
where $1\leqslant c_1<\dots<c_p\leqslant n$ and
$1\leqslant b_1<\dots<b_{m-2p}\leqslant 2n$ with the condition that
$b_i+b_j\ne 2n+1$ for all $i$ and $j$.
Applying the corresponding operator \eqref{qop}
to a basis vector
of this form we find that the coefficient
of this vector in the expansion equals
$r!{\displaystyle\binom{p}{r}}$.
Since $h_{i^{\tss\prime}}=-h_i$, the contributions of the
basis vectors with $m-2p>0$ in the respective sum
\eqref{sumis} add up to zero. This proves the
first part of the proposition.
Now suppose
that $m=2k$ is even. Then the coefficient of the monomial
$h_{i_1}^{2}\dots h_{i_k}^{2}$
with $1\leqslant i_1<\dots<i_k\leqslant n$
in the respective sum \eqref{sumis}
equals
\ben
(-1)^k\ts\sum_{r=0}^{k}
\binom{-n+2k-2}{r}^{-1}\binom{k}{r}=(-1)^k\ts\frac{-n+2k-1}{-n+k-1},
\een
completing the proof.
\epf

Now we will extend Proposition~\ref{prop:leadsksym}
to the values $n+1\leqslant m\leqslant 2n$.

\bpr\label{prop:recrels}
For $m\leqslant n$ suppose that $A,\dots,B,C$ are $m$
square matrices of size $2n$ with entries in $\CC$.
Then the following recurrence relations hold{\tss\rm :}
\beql{recid}
\tr\ts S^{(m)} A_1\dots B_{m-1}\tss C_m
=\frac{(n-m+1)(2n-m+3)}{m\tss(n-m+2)}\ts
\tr\ts S^{(m-1)} A_1\dots B_{m-1},
\eeq
if $C$ is the identity matrix, and
\beql{recskew}
\tr\ts S^{(m)} A_1\dots B_{m-1}\tss C_m=
-\frac{n-m+1}{m\tss(n-m+2)}\ts\sum_{i=1}^{m-1}\ts
\tr\ts S^{(m-1)} A_1\dots B_{m-1}\ts C_i,
\eeq
if $C$ satisfies the symmetry condition $C+C'=0$.
\epr

\bpf
Expression~\eqref{symjm} for the symmetrizer gives
the recurrence relation
\beql{reccsy}
S^{(m)}=\frac{1}{m\tss(\om+2m-4)}
\Big(1+\sum_{i=1}^{m-1}(s_{i\tss m}-\ep_{i\tss m})\Big)
\Big(\om+m-3+\sum_{i=1}^{m-1}(s_{i\tss m}-\ep_{i\tss m})\Big)
\tss S^{(m-1)}
\eeq
in the Brauer algebra $\Bc_m(\om)$.
Consider the corresponding operator
on the space \eqref{cntenpr} with $\om=-2n$ and
the action \eqref{braactsy}. Taking trace over the
$m$-th copy of $\End\CC^{2n}$ and using the relations
$\tr_m P_{i\tss m}=\tr_m Q_{i\tss m}=1$, $i<m$, we obtain
\ben
\tr_m\ts S^{(m)}
=\frac{(n-m+1)(2n-m+3)}{m\tss(n-m+2)}\ts S^{(m-1)}
\een
which gives \eqref{recid}. Similarly, if $C+C'=0$ then
\begin{multline}\non
\tr\ts S^{(m)} A_1\dots B_{m-1}\tss C_m\\
{}=\frac{1}{m\tss(-2n+2m-4)}\tss
\tr\ts S^{(m-1)}A_1\dots B_{m-1}\tss C_m
\Big({-1}+\sum_{i=1}^{m-1}\Phi_{i\tss m}\Big)
\Big(2n-m+3+\sum_{i=1}^{m-1}\Phi_{i\tss m}\Big),
\end{multline}
where we used the cyclic property of trace and set
$\Phi_{i\tss m}=P_{i\tss m}-Q_{i\tss m}$. For the partial trace we
calculate
\ben
\tr_m C_m\Big({-1}+\sum_{i=1}^{m-1}\Phi_{i\tss m}\Big)
\Big(2n-m+3+\sum_{i=1}^{m-1}\Phi_{i\tss m}\Big)=
2\tss(n-m+1)\ts\sum_{i=1}^{m-1}C_i
\een
thus proving \eqref{recskew}.
\epf

Proposition~\ref{prop:recrels} implies that if each of the $m$
matrices $A,\dots,C$ is either the identity matrix
or satisfies the symmetry condition, then for any
$1\leqslant l\leqslant\lfloor m/2\rfloor$ the expression
\beql{traa}
\frac{1}{n-m+1}\ts\tr\ts S^{(m)} A_1\dots  C_m
\eeq
can be written as a linear combination of traces of the form
\ben
\frac{1}{n-m+l+1}\ts\tr\ts S^{(m-l)} D_1\dots E_{m-l},
\een
where each of the matrices $D,\dots,E$
is a product of some of the matrices $A,\dots,C$.
Now observe that if $n+1\leqslant m\leqslant 2n$ then
the symmetrizer $S^{(m-l)}$ is well-defined for
$l=m-n$ so that the expression \eqref{traa} may be {\it defined\/}
as being equal to this linear combination.
To see that this evaluation of
\eqref{traa} is well-defined, it suffices
to note that various linear
combinations obtained by different applications of the
recurrence relations \eqref{recid}
and \eqref{recskew} agree for any fixed $m$ and infinitely
many values of $n$ with $n\geqslant m$.

With this interpretation of \eqref{traa},
relation \eqref{deftrsp} defines
the elements $\phi_{m\tss k}$ for all
values $1\leqslant m\leqslant 2n$. Moreover, we can now
extend Proposition~\ref{prop:leadsksym} to all these values
of the parameters.

\bco\label{cor:leadsp}
In the symplectic case with $1\leqslant m\leqslant 2n$,
the expression
\ben
\frac{1}{n-m+1}\ts\tr\ts S^{(m)} X_1\dots X_m
\een
is zero
if $m$ is odd, while for even values $m=2k$
we have
\ben
\frac{1}{n-2k+1}\ts\tr\ts S^{(2k)} X_1\dots X_{2k}
=\frac{(-1)^k}{n-k+1}
\sum_{1\leqslant i_1<\dots< i_k\leqslant n}
\ts h_{i_1}^2\dots h_{i_k}^2.
\een
\eco

\bpf
For a fixed value of $m$ regard $n$ as an integer parameter.
By Proposition~\ref{prop:leadsksym} the relations hold
for infinitely many values of $n$ with $n\geqslant m$, and so they
hold for
the values $m/2\leqslant n<m$ as well.
\epf

\section{Proof of the Main Theorem}
\label{sec:pmt}
\setcounter{equation}{0}

Recall that the dual Coxeter number for the Lie algebra
$\g_N$ is given by
\ben
h^{\vee}=\begin{cases}
N-2\quad&\text{for $\g_N=\oa_N$,}\\
n+1\quad&\text{for $\g_N=\spa_{2n}$.}
\end{cases}
\een
From now on we will assume that $\ka=-h^{\vee}$
so that $\ka$ is at
the {\it critical level\/}.

\subsection{Elements $\phi_{m\tss k}$}
\label{subsec:ephi}

To prove the first part of the theorem
(with the assumption $m\leqslant n$ in the symplectic case)
it will be sufficient to verify that
for all $i,j$
\beql{annih}
F_{ij}[0]\ts \tr\ts S^{(m)} \Phi_1\dots \Phi_m=
F_{ij}[1]\ts \tr\ts S^{(m)} \Phi_1\dots \Phi_m=0
\eeq
in the $\wh\g_N$-module
$V_{-h^{\vee}}(\g_N)\ot \CC[\tau]$.
Consider
the tensor product algebra
\beql{tenprext}
\End\CC^N\ot\dots\ot\End\CC^N\ot\U
\eeq
with $m+1$ copies of $\End\CC^N$ labeled by $0,1,\dots,m$,
where $\U$ stands for the universal enveloping algebra
$\U(\wh\g_N\oplus\CC\tau)$.
Relations
\eqref{annih} can now be written in the equivalent form
\beql{aanih}
F[0]_0\ts\tr\ts S^{(m)} \Phi_1\dots \Phi_m=0
\Fand F[1]_0\ts\tr\ts S^{(m)} \Phi_1\dots \Phi_m=0
\eeq
modulo the left ideal of $\U$ generated by
$\g_N[t]$ and $K+h^{\vee}$.
To verity the first relation, note that by \eqref{deffaff}
we have
\beql{fophi}
[F[0]_0,\Phi_i]=(P_{0\tss i}-Q_{0\tss i})\ts\Phi_i
-\Phi_i\ts(P_{0\tss i}-Q_{0\tss i}).
\eeq
Hence
\ben
\bal
F[0]_0\ts\tr\ts S^{(m)} \Phi_1\dots \Phi_m
{}=\sum_{i=1}^m \tr\ts S^{(m)} \Phi_1\dots
\Big((P_{0\tss i}-Q_{0\tss i})\ts\Phi_i
-\Phi_i\ts(P_{0\tss i}-Q_{0\tss i})\Big)\dots\Phi_m&\\
{}=\tr\ts S^{(m)} \Big(\sum_{i=1}^m(P_{0\tss i}-Q_{0\tss i})\Big)
\Phi_1\dots\Phi_m-\tr\ts S^{(m)}\ts
\Phi_1\dots\Phi_m\ts \Big(\sum_{i=1}^m(P_{0\tss i}-Q_{0\tss i})\Big).&
\eal
\een
The sum $\sum_{i=1}^m(P_{0\tss i}-Q_{0\tss i})$
commutes with the action of
the symmetrizer $S^{(m)}$ so that the expression
is equal to zero by the cyclic property of trace.

For the proof of the second relation in
\eqref{aanih} use the following consequence
of \eqref{deffaff}:
\ben
[F[1]_0,\Phi_i]=F[0]_0
+(P_{0\tss i}-Q_{0\tss i})\ts F[0]_i-F[0]_i\ts(P_{0\tss i}-Q_{0\tss i})
+(P_{0\tss i}-Q_{0\tss i})\ts K',
\een
where $K'=K$ in the orthogonal case and
$K'=2\tss K$ in the symplectic case. We have
\ben
\bal
F[1]_0&\ts\tr\ts S^{(m)} \Phi_1\dots \Phi_m
=\sum_{i=1}^m \tr\ts S^{(m)} \Phi_1\dots\Phi_{i-1}\\
{}&\times\Big(F[0]_0
+(P_{0\tss i}-Q_{0\tss i})\ts F[0]_i-F[0]_i\ts(P_{0\tss i}-Q_{0\tss i})
+(P_{0\tss i}-Q_{0\tss i})\ts K'\Big)\dots\Phi_m
\eal
\een
and applying \eqref{fophi} we can write this expression as
\ben
\bal
&\sum_{1\leqslant i<j\leqslant m} \tr\ts S^{(m)} \Phi_1\dots
\wh \Phi_i\dots \Big((P_{0\tss j}-Q_{0\tss j})\ts\Phi_j
-\Phi_j\ts(P_{0\tss j}-Q_{0\tss j})\Big)\dots\Phi_m\\
&{}+\sum_{1\leqslant i<j\leqslant m} \tr\ts S^{(m)}
(P_{0\tss i}-Q_{0\tss i})\ts\Phi_1\dots
\wh \Phi_i\dots \Big((P_{ij}-Q_{ij})\ts\Phi_j
-\Phi_j\ts(P_{ij}-Q_{ij})\Big)\dots\Phi_m\\
&{}-\sum_{1\leqslant i<j\leqslant m} \tr\ts S^{(m)}
\Phi_1\dots
\wh \Phi_i\dots \Big((P_{ij}-Q_{ij})\ts\Phi_j
-\Phi_j\ts(P_{ij}-Q_{ij})\Big)\dots\Phi_m\ts(P_{0\tss i}-Q_{0\tss i})\\
&{}+K'\ts\sum_{i=1}^m\ts \tr\ts S^{(m)}
(P_{0\tss i}-Q_{0\tss i})\ts\Phi_1\dots
\wh \Phi_i\dots \Phi_m,
\eal
\een
where hats indicate that the corresponding symbols
should be skipped. Now we simplify it by using
the cyclic property of trace,
conjugations by elements of $\Sym_m$ and relations \eqref{bradr}.
In the orthogonal case, for the first sum we have
\ben
\bal
\sum_{1\leqslant i<j\leqslant m} \tr\ts S^{(m)} \Phi_1\dots
&\wh \Phi_i\dots (P_{0\tss j}-Q_{0\tss j})\ts\Phi_j
\dots\Phi_m\\
{}&=\sum_{1\leqslant i<j\leqslant m} \tr\ts S^{(m)}
P_{m-1,m}\dots P_{i,i+1}\ts
\Phi_1\dots
\wh \Phi_i\dots (P_{0\tss j}-Q_{0\tss j})\ts\Phi_j
\dots\Phi_m\\
{}&=\sum_{1\leqslant i<j\leqslant m} \tr\ts S^{(m)}
(P_{0\ts j-1}-Q_{0\ts j-1})\ts
\Phi_1\dots
\Phi_{m-1}\ts P_{m-1,m}\dots P_{i,i+1}\\
{}&=\sum_{i=1}^{m-1} i\ts\tr\ts  S^{(m)}
(P_{0\tss i}-Q_{0\tss i})\ts
\Phi_1\dots
\Phi_{m-1}.
\eal
\een
In the symplectic case,
$S^{(m)}=-S^{(m)}\ts P_{j\ts j+1}$ so that
the final expression in this chain of equalities will have
the same form. In both cases, applying similar transformations
to all remaining sums we will get
\ben
\bal
F[1]_0\ts\tr\ts S^{(m)} \Phi_1\dots \Phi_m
{}&=\sum_{i=1}^{m-1} i\ts\tr\ts  S^{(m)}
(P_{0\tss i}-Q_{0\tss i})\ts
\Phi_1\dots
\Phi_{m-1}\\
{}&-\sum_{i=1}^{m-1} i\ts\tr\ts(P_{0\tss i}-Q_{0\tss i})\ts  S^{(m)}
\Phi_1\dots
\Phi_{m-1}\\
&+\sum_{i=1}^{m-1} i\ts\tr\ts S^{(m)}
(P_{0\tss m}-Q_{0\tss m})
(P_{i\tss m}-Q_{i\tss m})\ts
\Phi_1\dots
\Phi_{m-1}\\
&+\sum_{i=1}^{m-1} i\ts\tr
\ts(P_{i\tss m}-Q_{i\tss m})\ts(P_{0\tss m}-Q_{0\tss m})
\ts S^{(m)}
\Phi_1\dots
\Phi_{m-1}\\
&+\big(m\tss K'\mp m(m-1)\big)\ts\tr
\ts S^{(m)}(P_{0\tss m}-Q_{0\tss m})\ts
\Phi_1\dots
\Phi_{m-1},
\eal
\een
where the upper sign in the double sign corresponds to
the orthogonal case and lower sign to the symplectic case.
Furthermore, observe that
\ben
S^{(m)}
(P_{0\tss m}-Q_{0\tss m})\ts P_{i\tss m}
=S^{(m)} P_{i\tss m}
\ts(P_{0\tss i}-Q_{0\tss i})={}\pm S^{(m)}
(P_{0\tss i}-Q_{0\tss i}),
\een
and that
\beql{spqze}
S^{(m)}
(P_{0\tss m}-Q_{0\tss m})\ts Q_{i\tss m}=0
\eeq
since
\ben
\bal
S^{(m)}
(P_{0\tss m}-Q_{0\tss m})\ts Q_{i\tss m}&={}\pm S^{(m)}
(P_{0\tss m}-Q_{0\tss m})\ts P_{i\tss m}\ts Q_{i\tss m}\\
&=S^{(m)}
(P_{0\tss i}-Q_{0\tss i})\ts Q_{i\tss m}=
-S^{(m)}
(P_{0\tss m}-Q_{0\tss m})\ts Q_{i\tss m}.
\eal
\een
Hence, simplifying the expressions further we get
\ben
\bal
F[1]_0\ts\tr\ts S^{(m)} \Phi_1\dots \Phi_m
{}&=\sum_{i=1}^{m-1} 2\ts i\ts\tr\ts  S^{(m)}
(P_{0\tss i}-Q_{0\tss i})\ts
\Phi_1\dots
\Phi_{m-1}\\
&+\big(m\tss K'-m(m-1)\big)\ts\tr
\ts S^{(m)}(P_{0\tss m}-Q_{0\tss m})\ts
\Phi_1\dots
\Phi_{m-1}
\eal
\een
in the orthogonal case, and
\ben
\bal
F[1]_0\ts\tr\ts S^{(m)} \Phi_1\dots \Phi_m
{}&=-\sum_{i=1}^{m-1} 2\ts i\ts
\tr\ts(P_{0\tss i}-Q_{0\tss i})\ts  S^{(m)}
\Phi_1\dots
\Phi_{m-1}\\
&+\big(m\tss K'+m(m-1)\big)\ts\tr
\ts S^{(m)}(P_{0\tss m}-Q_{0\tss m})\ts
\Phi_1\dots
\Phi_{m-1}
\eal
\een
in the symplectic case.
As a next step, calculate the partial trace $\tr_m$ on the right
hand sides with respect to the $m$-th copy of
$\End\CC^N$ in \eqref{tenprext} with the use of the following
lemma.

\ble\label{lem:trsym}
We have
\ben
\tr_m\ts S^{(m)}={}\pm\frac{(\om+m-3)(\om+2m-2)}{m\ts(\om+2m-4)}\ts
S^{(m-1)}
\een
and
\ben
\tr_m\ts S^{(m)}(P_{0\tss m}-Q_{0\tss m})=
{}\pm\frac{\om+2m-2}{m\ts(\om+2m-4)}\ts
S^{(m-1)}\sum_{i=1}^{m-1}(P_{0\tss i}-Q_{0\tss i}),
\een
where $\om$ equals $N$ and $-N$ in the orthogonal and symplectic
case, respectively.
\ele

\bpf
This is derived from the recurrence relation \eqref{reccsy}
in the same way as in the proof of Proposition~\ref{prop:recrels}.
\epf

Applying Lemma~\ref{lem:trsym} we come to the relations
\begin{multline}\non
F[1]_0\ts\tr\ts S^{(m)} \Phi_1\dots \Phi_m=
\frac{N+2\tss m-2}{m\ts(N+2\tss m-4)}\\
{}{}\times\Big((N+m-3)
\sum_{i=1}^{m-1} 2\ts i\ts\tr\ts  S^{(m-1)}
(P_{0\tss i}-Q_{0\tss i})\ts
\Phi_1\dots
\Phi_{m-1}\\
+m\tss \big(K'-m+1\big)
\sum_{i=1}^{m-1}\tr
\ts S^{(m-1)}(P_{0\tss i}-Q_{0\tss i})\ts
\Phi_1\dots
\Phi_{m-1}\Big)
\end{multline}
in the orthogonal case, and
\begin{multline}\non
F[1]_0\ts\tr\ts S^{(m)} \Phi_1\dots \Phi_m=
-\frac{N-2\tss m+2}{m\ts(N-2\tss m+4)}\\
{}{}\times\Big((N-m+3)
\sum_{i=1}^{m-1} 2\ts i\ts\tr\ts
(P_{0\tss i}-Q_{0\tss i})\ts S^{(m-1)}
\Phi_1\dots
\Phi_{m-1}\\
+m\tss \big(K'+m-1\big)
\sum_{i=1}^{m-1}\tr\ts
(P_{0\tss i}-Q_{0\tss i})\ts S^{(m-1)}
\Phi_1\dots
\Phi_{m-1}\Big)
\end{multline}
in the symplectic case, where the traces on the right hand sides
are now taken over the $m-1$ copies of $\End\CC^N$.

\ble\label{lem:treq}
For any $1\leqslant i<j\leqslant m$ we have
\ben
\tr\ts S^{(m)}
(P_{0\tss i}-Q_{0\tss i})\ts
\Phi_1\dots\Phi_{m}=\tr\ts S^{(m)}
(P_{0\tss j}-Q_{0\tss j})\ts
\Phi_1\dots\Phi_{m}
\een
in the orthogonal case, and
\ben
\tr\ts (P_{0\tss i}-Q_{0\tss i})\ts S^{(m)}
\Phi_1\dots\Phi_{m}=\tr
\ts(P_{0\tss j}-Q_{0\tss j})\ts S^{(m)}
\Phi_1\dots\Phi_{m}
\een
in the symplectic case, where the traces are taken over
all $m$ copies of $\End\CC^N$.
\ele

\bpf
We may assume that $m=2$. In the orthogonal case we have
\ben
\tr\ts S^{(2)}\ts (P_{0\tss 1}-Q_{0\tss 1}-P_{0\tss 2}+Q_{0\tss 2})
\ts\Phi_1\Phi_2
={}-\tr\ts S^{(2)}\ts (P_{0\tss 1}-Q_{0\tss 1}-P_{0\tss 2}+Q_{0\tss 2})
\ts\Phi_2\Phi_1,
\een
and by \eqref{deffaff}
\begin{multline}\non
\tr\ts S^{(2)}\ts (P_{0\tss 1}-Q_{0\tss 1}-P_{0\tss 2}+Q_{0\tss 2})
\ts[\Phi_1,\Phi_2]
=\tr\ts S^{(2)}\ts (P_{0\tss 1}-Q_{0\tss 1}-P_{0\tss 2}+Q_{0\tss 2})\\
{}\times
\Big((F[-2]_1-F[-2]_2)(P_{12}-1)
+Q_{12}\ts F[-2]_1-F[-2]_1\ts Q_{12}\Big)
\end{multline}
which is zero by \eqref{bradr}, \eqref{spqze}
and the cyclic property of trace.
The calculation in the symplectic case is similar;
it relies on the fact that the trace
\begin{multline}\non
\tr\ts S^{(2)}[\Phi_1,\Phi_2]\ts
(P_{0\tss 1}-Q_{0\tss 1}-P_{0\tss 2}+Q_{0\tss 2})\\
=\tr\ts S^{(2)}\Big((P_{12}+1)(F[-2]_2-F[-2]_1)
+Q_{12}\ts F[-2]_1-F[-2]_1\ts Q_{12}\Big)
\ts (P_{0\tss 1}-Q_{0\tss 1}-P_{0\tss 2}+Q_{0\tss 2})
\end{multline}
is zero.
\epf

By Lemma~\ref{lem:treq}, up to a numerical factor,
$F[1]_0\ts\tr\ts S^{(m)} \Phi_1\dots \Phi_m$ equals
\ben
(K'+N\mp2)\ts \sum_{i=1}^{m-1} \ts\tr\ts
 S^{(m-1)}(P_{0\tss i}-Q_{0\tss i})\ts
\Phi_1\dots \Phi_{m-1}.
\een
At the critical level, $K'+N-2=K+N-2=0$
and $K'+N+2=2\tss(K+n+1)=0$ in the orthogonal and symplectic
case, respectively. This completes the proof
of the first part of the Main Theorem
in the orthogonal case, while in the symplectic case
the proof is complete under
the assumption $m\leqslant n$.

Now we continue the argument in the symplectic case.
Recalling that $\Phi=\tau+F[-1]$, we can write
$\Phi_1\dots \Phi_m$ as a polynomial in $\tau$
whose coefficients are linear combinations of
products of the form $F[r_1]_{p_1}\dots F[r_s]_{p_s}$
with the conditions
$1\leqslant p_1<\dots<p_s\leqslant m$ and $r_i=-1,-2,\dots$.
Using the cyclic property of trace and
applying conjugations by elements of $\Sym_m$ we find
that for $n\geqslant m$ the element
$\tr\ts S^{(m)} F[r_1]_{p_1}\dots F[r_s]_{p_s}$
coincides with $\tr\ts S^{(m)} F[r_1]_1\dots F[r_s]_s$.
Since the matrix $F[r]$ satisfies
$F[r]+F[r]^{\tss\prime}=0$,
it can be written in the form
\ben
F[r]=\sum_{i,j=1}^{2n} e_{ij}\ot F_{ij}[r]
=\frac12\ts\sum_{i,j=1}^{2n} e_{ij}\ot
(F^{}_{ij}[r]-F_{ij}[r]^{\tss\prime})
=\frac12\ts\sum_{i,j=1}^{2n} f_{ij}\ot F_{ij}[r],
\een
where $f_{ij}$ is defined
by
\ben
f_{ij}=e_{ij}-\sum_{k,\tss l=1}^{2n}
\overline g^{}_{ik}\tss  g^{}_{lj}\tss e_{lk}\in\End\CC^{2n}.
\een
Hence
\begin{multline}\non
\frac{1}{n-m+1}\ts\tr\ts S^{(m)} F[r_1]_1\dots F[r_s]_s=\\
\sum_{i_1,j_1,\dots,i_s,j_s}
\Big(\frac{1}{2^s(n-m+1)}\ts\tr\ts S^{(m)}f_{i_1j_1}\ot\dots
\ot f_{i_sj_s}\ot 1\ot\dots
\ot 1\Big)\ot
F_{i_1j_1}[r_1]\dots F_{i_sj_s}[r_s].
\end{multline}
Note that the expression in the brackets has the form \eqref{traa}
so it is well-defined for the values of $m$ with
$n+1\leqslant m\leqslant 2n$. To show that the coefficients
$\phi_{m\tss k}$ defined in \eqref{deftrsp} are
Segal--Sugawara
vectors for $\spa_{2n}$, note that for
each fixed value of $m$
the relations
\eqref{aanih} hold for infinitely many
values $n\geqslant m$ of the parameter $n$. This implies that
they hold for all values $m/2\leqslant n<m$.

To prove the remaining parts of the Main Theorem,
note that the symbols of
the elements $\phi_{2k\ts 2k}$ with $k=1,\dots,n$
defined in \eqref{deftr} and \eqref{deftrsp} respectively coincide with
the images of the leading terms of the invariants
$\tr\ts S^{(2k)} F_1\dots F_{2k}\in\U(\g_N)^{\Gr_N}$
introduced in Sec.~\ref{subsec:cm}
(with the additional factor $(n-2k+1)^{-1}$
in the symplectic case),
under the embedding $F_{ij}\mapsto F_{ij}[-1]$.
Due to Proposition~\ref{prop:leadsym} and
Corollary~\ref{cor:leadsp},
the images of the leading terms in the algebra of
symmetric polynomials in $h_1^2,\dots,h_n^2$ coincide,
up to numerical factors, with the complete and elementary symmetric
functions, respectively. The theorem is now proved
for $\g_N$ of types $B$ and $C$.

\subsection{Pfaffians}
\label{subsec:pf}

To complete the argument
in the case $\g_N=\oa_{2n}$ recall
that the element $\phi_n'=\Pf \wt F[-1]$ in type $D$
is defined in \eqref{genpf}.

\ble\label{lem:pfaff}
The element $\phi'_n=\Pf \wt F[-1]$ is a Segal--Sugawara vector for $\wh\oa_{2n}$.
\ele

\bpf
Suppose first that
the bilinear form defined in \eqref{formg} corresponds
to the identity matrix $G=1$. The generators $F_{ij}$ of $\oa_{2n}$
are given by \eqref{fijelem}.
In this case $\wt F[-1]=F[-1]$, and the commutation relations
of the Lie algebra $\wh\oa_{2n}$ simplify to
\ben
\bal
\big[F_{ij}[r],F_{kl}[s]\big]&=\de_{kj}\ts F_{il}[r+s]-\de_{il}\ts F_{kj}[r+s]
-\de_{ki}\ts F_{jl}[r+s]+\de_{jl}\ts F_{ki}[r+s]\\
{}&+r\tss\de_{r,-s}\big(\de_{kj}\ts\de_{il}-\de_{ki}\ts\de_{jl}\big)\ts K.
\eal
\een
In particular,
the elements $F_{ij}[-1]$ and $F_{kl}[-1]$ of the Lie algebra $\wh\oa_{2n}$ commute
if the indices $i,j,k,l$ are distinct. Therefore, we can write the formula
for the Pfaffian in the form
\beql{pf}
\Pf F[-1]=\sum_{\si}\sgn\si\cdot
F_{\si(1)\ts\si(2)}[-1]\dots F_{\si(2n-1)\ts\si(2n)}[-1],
\eeq
summed over the elements $\si$ of the subset $\Ac_{2n}\subset \Sym_{2n}$
which consists of the permutations with the properties
$\si(2k-1)<\si(2k)$ for all $k=1,\dots,n$ and $\si(1)<\si(3)<\dots<\si(2n-1)$.

We need to show that $\oa_{2n}[t]\ts \phi'_n=0$ in the vacuum module $V_{-h^{\vee}}(\oa_{2n})$.
It suffices to verify that for all $i,j$,
\beql{annihpf}
F_{ij}[0]\ts \Pf F[-1]=F_{ij}[1]\ts \Pf F[-1]=0.
\eeq
Note that for any permutation $\pi\in\Sym_{2n}$ the mapping
\ben
F_{ij}[r]\mapsto F_{\pi(i)\ts\pi(j)}[r],\qquad K\mapsto K
\een
defines an automorphism of the Lie algebra $\wh\oa_{2n}$. Moreover,
the image of $\Pf F[-1]$ under its extension to $\U(\wh\oa_{2n})$ coincides
with $\sgn\pi\cdot\Pf F[-1]$. Hence, it is enough to verify \eqref{annihpf}
for $i=1$ and $j=2$.
Observe that $F_{12}[0]$ commutes with all summands in \eqref{pf} with
$\si(1)=1$ and $\si(2)=2$. Suppose now that $\si\in\Ac_{2n}$
is such that $\si(2)>2$. Then $\si(3)=2$ and $\si(4)>2$.
In $V_{-h^{\vee}}(\oa_{2n})$ we have
\ben
\bal
F_{12}[0]\ts F_{1\ts \si(2)}[-1]\tss F_{2\ts \si(4)}[-1]&\dots F_{\si(2n-1)\ts\si(2n)}[-1]\\
{}={}&{}-F_{2\ts \si(2)}[-1]\tss F_{2\ts \si(4)}[-1]\dots F_{\si(2n-1)\ts\si(2n)}[-1]\\
{}&{}+F_{1\ts \si(2)}[-1]\tss F_{1\ts \si(4)}[-1]\dots F_{\si(2n-1)\ts\si(2n)}[-1].
\eal
\een
Set $i=\si(2)$ and $j=\si(4)$. Note that the permutation $\si'=\si\ts (2\ts 4)$ also belongs
to the subset $\Ac_{2n}$, and $\sgn\si'=-\sgn\si$. We have
\begin{multline}
-F_{2\ts i}[-1]\tss F_{2\ts j}[-1]+F_{1\ts i}[-1]\tss F_{1\ts j}[-1]
+F_{2\ts j}[-1]\tss F_{2\ts i}[-1]-F_{1\ts j}[-1]\tss F_{1\ts i}[-1]\\
{}=F_{ij}[-2]-F_{ij}[-2]=0.
\non
\end{multline}
This implies that the terms in the expansion of $F_{12}[0]\ts\Pf F[-1]$
corresponding to pairs of the form $(\si,\si')$ cancel pairwise. Thus,
$F_{12}[0]\ts\Pf F[-1]=0$.
Now we verify that
\beql{fone}
F_{12}[1]\ts\Pf F[-1]=0.
\eeq
Consider first the summands in \eqref{pf} with
$\si(1)=1$ and $\si(2)=2$. In $V_{-h^{\vee}}(\oa_{2n})$ we have
\begin{align}
F_{12}[1]\ts F_{1\ts 2}[-1]\tss F_{\si(3)\ts \si(4)}[-1]&
\dots F_{\si(2n-1)\ts\si(2n)}[-1]
\non\\
{}&=-K\ts F_{\si(3)\ts \si(4)}[-1]\dots F_{\si(2n-1)\ts\si(2n)}[-1].
\label{fsi}
\end{align}
Furthermore, let $\tau\in\Ac_{2n}$
with $\tau(2)>2$. Then $\tau(3)=2$ and $\tau(4)>2$. We have
\begin{align}
F_{12}[1]\ts F_{1\ts \tau(2)}[-1]\tss F_{2\ts \tau(4)}[-1]&
\dots F_{\tau(2n-1)\ts\tau(2n)}[-1]
\non\\
{}&=-F_{2\ts \tau(2)}[0]\tss F_{2\ts \tau(4)}[-1]\dots F_{\tau(2n-1)\ts\tau(2n)}[-1]
\non\\
{}&= F_{\tau(2)\ts \tau(4)}[-1]\dots F_{\tau(2n-1)\ts\tau(2n)}[-1].
\label{taumonom}
\end{align}
Note that for any given $\si$ occurring in \eqref{fsi} there are exactly
$2n-2$ elements $\tau$ such that the product \eqref{taumonom}
coincides with $F_{\si(3)\ts \si(4)}[-1]\dots F_{\si(2n-1)\ts\si(2n)}[-1]$,
up to a sign. Hence, taking the signs
of permutations into account, we find that
\ben
F_{12}[1]\ts \Pf F[-1]
=(-K-2n+2)\ts \sum_{\si}\ts \sgn\si\cdot
F_{\si(3)\ts \si(4)}[-1]\dots F_{\si(2n-1)\ts\si(2n)}[-1],
\een
summed over $\si\in\Ac_{2n}$ with $\si(1)=1$ and $\si(2)=2$.
Since $-K-2n+2=0$ at the critical level, we get \eqref{fone}.

Finally, suppose that the bilinear form \eqref{formg} corresponds to
an arbitrary nondegenerate symmetric matrix $G$. Fix a matrix $A$ such that
$AA^t=G$. To distinguish the presentation of the Lie algebra $\oa_{2n}$
associated with the identity matrix as above, we
denote by $F^{\circ}_{ij}$ the generators of $\oa_{2n}=\oa_{2n}(1)$
used in that presentation.
An isomorphism between the presentations $\oa_{2n}(G)$ and
$\oa_{2n}(1)$ of the Lie algebra $\oa_{2n}$, associated with the respective matrices
$G$ and $1$ can be given by
$F\mapsto AF^{\circ}A^{-1}$.
This isomorphism extends naturally to the corresponding presentations of the affine
Kac--Moody algebra $\wh\oa_{2n}$ so that $\wt F[r]\mapsto AF^{\circ}[r]A^t$
for any $r\in\ZZ$, where $\wt F[r]=F[r]\tss G$.
Since the Pfaffian $\Pf \wt F[-1]$ is found by the
expansion
\ben
\frac{\Fc^{\tss n}}{n\tss!}=
\big(\tss e_{1}\wedge\dots\wedge e_{2n}\big)\ot\Pf \wt F[-1],
\een
with
\ben
\Fc=\sum_{i<j}\ (e_i\wedge e_j)\ot \wt F_{ij}[-1]\in\La(\CC^{2n})\ot\U(\wh\oa_{2n}),
\een
we conclude that the image of $\Pf \wt F[-1]$ under the isomorphism
between the presentations of $\wh\oa_{2n}$ is found by
\ben
\Pf \wt F[-1]\mapsto \det A\cdot \Pf F^{\circ}[-1].
\een
This implies that $\Pf \wt F[-1]$ is a Segal--Sugawara vector
for the presentation of $\wh\oa_{2n}$ associated with an arbitrary
nondegenerate symmetric matrix $G$.
\epf

Clearly, the leading term of the Segal--Sugawara vector $\phi'_n$
coincides with the image of the Pfaffian invariant
in the symmetric algebra $\Sr(\oa_{2n})$ under the embedding
$F_{ij}\mapsto F_{ij}[-1]$. This completes the proof of the Main Theorem
in the case $D$.

\subsection{Corollaries: traces over $\sll_2$-singular vectors}
\label{subsec:cor}

Now we use the results of Sec.~\ref{subsec:sv}
to derive expressions for the Segal--Sugawara vectors
$\phi^{}_{m\tss k}$ based on the use of the Howe dual
pairs $(\sll_2,\g_N)$. In the orthogonal case,
using the notation as in \eqref{tenprext}, for each $m$
introduce the element
$\Phi^{(m)}\in\End\Pc_N^{\tss m}\ot\U$
by setting
\ben
\Phi^{(m)}:z_{j_1}\dots z_{j_m} \mapsto
\sum_{i_1\leqslant\dots\leqslant i_m}
\ts z_{i_1}\dots z_{i_m}\ot
\Phi^{\tss i_1,\dots,i_m}_{\tss j_1,\dots,j_m}
\een
where
\ben
\Phi^{\tss i_1,\dots,i_m}_{\tss j_1,\dots,j_m}
=\frac{1}{\al_1!\dots\al_N!\ts m!}\sum_{\si,\pi\in\Sym_m}
\ts \Phi_{i_{\si(1)}j_{\pi(1)}}\dots \Phi_{i_{\si(m)}j_{\pi(m)}}
\een
and $\al_i$ is the
multiplicity of $i$ in the multiset
$\{i_1,\dots,i_m\}$. The Main Theorem
and Proposition~\ref{prop:expsm} imply the
following corollary.

\bco\label{cor:singv}
The Segal--Sugawara vectors
$\phi^{}_{m\tss k}$ can be found from
the expansion
\ben
\tr\ts p\ts\Phi^{(m)}=\phi^{}_{m\tss0}\ts
\tau^m+\phi^{}_{m\tss1}\ts\tau^{m-1}
+\dots+\phi^{}_{m\tss m}
\een
with the trace taken over the subspace of $\sll_2$-singular
vectors in $\Pc_N^{\tss m}$.
\qed
\eco

Similarly, in the symplectic case introduce
the element $\Phi^{(m)}\in\End\La_{2n}^m\ot\U$ by
\ben
\Phi^{(m)}:\ze_{j_1}\wedge\dots\wedge\ze_{j_m}
\mapsto \sum_{i_1<\dots<i_m}\ts
\ze_{i_1}\wedge\dots\wedge\ze_{i_m}\ot
\Phi^{\tss i_1,\dots,i_m}_{\tss j_1,\dots,j_m},
\een
where
\ben
\Phi^{\tss i_1,\dots,i_m}_{\tss j_1,\dots,j_m}
=\frac{1}{m!}\ts\sum_{\si,\pi\in\Sym_m}\sgn\si\pi
\ts \Phi_{i_{\si(1)}j_{\pi(1)}}\dots \Phi_{i_{\si(m)}j_{\pi(m)}}.
\een
Applying the Main Theorem
and Proposition~\ref{prop:expsmsp} we get
the following corollary.

\bco\label{cor:singvsp}
The Segal--Sugawara vectors
$\phi^{}_{m\tss k}$ with $m\leqslant n$ can be found from
the expansion
\ben
\frac{1}{n-m+1}\ts\tr\ts p\ts\Phi^{(m)}=\phi^{}_{m\tss0}\ts
\tau^m+\phi^{}_{m\tss1}\ts\tau^{m-1}
+\dots+\phi^{}_{m\tss m}
\een
with the trace taken over the subspace of $\sll_2$-singular
vectors in $\La_{2n}^{\tss m}$. Moreover, for any fixed value
of $m$ the vectors $\phi^{}_{m\tss k}$ are well-defined
for the specializations of $n$ to all values $m/2\leqslant n<m$.
\qed
\eco

\section{Gaudin Hamiltonians and Bethe subalgebras}
\label{sec:hgh}
\setcounter{equation}{0}

Apply the state-field correspondence map $Y$
to the Segal--Sugawara vectors
provided by the Main Theorem.
Introduce the matrix
$F(z)=[F_{ij}(z)]$, where
the fields $F_{ij}(z)$ are defined in \eqref{basfi}.
Then in the orthogonal case
all Fourier coefficients of the fields
$f_{m\tss i}(z)$ defined by the decompositions
of the normally ordered trace
\beql{expno}
:\tr\ts S^{(m)} \big(\di_z+F_1(z)\big)\dots \big(\di_z+F_m(z)\big):{}=
f_{m\tss 0}(z)\ts\di_z^{\tss m}+f_{m\tss 1}(z)\ts\di_z^{\tss m-1}
+\dots+f_{m\tss m}(z)
\eeq
belong to the center
of the local completion of
the universal enveloping algebra $\U(\wh\g_N)$;
that is, they are
{\it Sugawara operators\/} for $\wh\g_N$; see \cite{ff:ak}.
The same holds in the symplectic case, where the left hand side
of \eqref{expno} should get the factor $(n-m+1)^{-1}$
and the values of $m$ restricted to $1\leqslant m\leqslant 2n$.
By the vacuum axiom of a vertex algebra, the application
of the fields $f_{m\tss i}(z)$
to the vacuum vector yields formal power series in $z$
with coefficients in
$\U(t^{-1}\g_N[t^{-1}])$. All coefficients of these
formal power series belong
to the center $\z(\wh\g_N)$ of the vertex algebra
$V_{-h^{\vee}}(\g_N)$. Moreover,
by the Main Theorem,
the center is generated
by these coefficients
(together with the additional series corresponding to
the Segal--Sugawara vector $\phi_n'$ in type $D$); see \cite{ff:ak},
\cite[Ch.~3]{f:lc}.
Thus we get an alternative presentation of the
generators of this commutative subalgebra of
$\U(t^{-1}\g_N[t^{-1}])$; see also a general result in \cite{r:uh}
which implies that this subalgebra is maximal commutative.
Explicitly, the coefficients of these formal power series are given
by the same expansions as in \eqref{expno}
by omitting the normal ordering signs and by
replacing $F(z)$ with the matrix
$F(z)_+=[F_{ij}(z)_+]$, where we use the standard notation
\ben
F_{ij}(z)_+=
\sum_{r=0}^{\infty} F_{ij}[-r-1]\tss z^{r},\qquad
F_{ij}(z)_-=
\sum_{r=0}^{\infty} F_{ij}[r]\tss z^{-r-1}.
\een
The same argument as in \cite[Sec.~3.2]{mr:mm} then yields
a corresponding family of commuting elements in $\U(\g_N[t])$; cf.~\cite[Sec.~4]{r:si}.
Indeed, given any $N\times N$ matrix $B=[b_{ij}]$ over $\CC$
with the condition $B+B^{\tss\prime}=0$, the commutation relations
between the series $F_{ij}(z)_-$
remain valid after the replacement
$
F_{ij}(z)_-\mapsto b_{ij}+F_{ij}(z)_-.
$
Introduce the matrix $F(z)_-=[F_{ij}(z)_-]$
and set $L(z)=\di_z-B-F(z)_-$.

\bco\label{cor:commsub}
The coefficients of all series
$l_{m\tss i}(z)$ with $m=1,2,\dots$
defined by the decompositions
\beql{composi}
\tr\ts S^{(m)} L_1(z)\dots L_m(z){}=
l_{m\tss 0}(z)\ts\di_z^{\tss m}+
l_{m\tss 1}(z)\ts\di_z^{\tss m-1}
+\dots+l_{m\tss m}(z),
\eeq
where in the symplectic case $m\leqslant 2n$
and the left hand side gets the factor $(n-m+1)^{-1}$,
generate a commutative subalgebra of $\U(\g_N[t])$.
\qed
\eco

Suppose now that $M$ is a finite-dimensional $\g_N$-module
and $a\in\CC$. Define
the corresponding evaluation $\g_N[t]$-module $M_a$
via the evaluation homomorphism
\beql{evalhom}
\ev_a:\g_N[t]\to \g_N,\qquad F_{ij}[r]\mapsto F_{ij}\ts a^r,
\eeq
that is,
$
F_{ij}(z)_-\mapsto F_{ij}/(z-a).
$
Given
finite-dimensional $\g_N$-modules $M^{(1)},\dots,M^{(p)}$
and complex numbers $a_1,\dots,a_p$,
the tensor product
$
M^{(1)}_{a_1}\ot\dots\ot M^{(p)}_{a_p}
$
becomes a $\g_N[t]$-module such that
the images of the entries of the matrix
$L(z)$ are found by
\ben
\ell_{ij}(z)=\de_{ij}\di_z-b_{ij}
-\sum_{s=1}^p\frac{F_{ij}^{(s)}}{z-a_s},
\een
where $F_{ij}^{(s)}$ denotes the image of $F_{ij}$
in $M^{(s)}$. Replacing
$L(z)$ by the matrix $[\ell_{ij}(z)]$
in \eqref{composi}
we obtain a family of commuting operators
in the tensor product module, thus producing higher
Gaudin Hamiltonians associated with $\g_N$;
cf. \cite{ct:qs},
\cite{ffr:gm}, \cite{fft:gm}, \cite{mtv:be}, \cite{r:si}.

Taking the evaluation homomorphism \eqref{evalhom}
with $a=0$ and applying the general
results on the ``shift of argument subalgebras" (see
\cite{fft:gm}, \cite{r:si} and \cite{t:ms}), we come to another corollary
providing explicit algebraically independent
generators of those subalgebras; cf. \cite{no:bs}. Write the decompositions
\eqref{composi} with the matrix
$L(z)=\di_z-B-F\tss z^{-1}$.
The components $l_{m\tss i}(z)$ defined by \eqref{composi}
are now polynomials in $z^{-1}$
with coefficients in $\U(\g_N)$.
In particular, let
\ben
l_{m\tss m}(z)=l^{\tss(0)}_{m\tss m}+l^{\tss(1)}_{m\tss m}\tss z^{-1}+\cdots
+ l^{\tss(m)}_{m\tss m}\tss z^{-m},\qquad l^{\tss(s)}_{m\tss m}\in\U(\g_N).
\een
Fix a Cartan subalgebra $\h_N$ of the Lie algebra $\g_N$. We will regard
any matrix $B$ with the property $B+B^{\tss\prime}=0$ as an element of $\g^*_N$.
If $\g_N=\oa_{2n}$, then
the matrix $\wt B=B\tss G$ is skew-symmetric and we define
the elements $p^{(r)}\in\U(\oa_{2n})$ by the expansion of the Pfaffian,
\beql{pfexp}
\Pf(\wt B+\wt F\tss z^{-1})=p^{(0)}+p^{(1)} z^{-1}+\cdots+p^{(n)} z^{-n},\qquad
\wt F=F\tss G.
\eeq

\bco\label{cor:commuea} {\rm(i)}\quad Suppose that
$\g_N=\oa_{2n+1}$ or $\g_N=\spa_{2n}$.
Then the coefficients of all polynomials
$l_{m\tss i}(z)$
generate a commutative subalgebra of $\U(\g_N)$.
If $B$ is a regular element of $\h^*_N$, then
this subalgebra is maximal commutative.
It is freely generated by the elements $l^{\tss(1)}_{2k\ts 2k},
\dots,l^{\tss(2k)}_{2k\ts 2k}$ with $k=1,\dots,n$.
\par
{\rm(ii)}\quad Suppose that
$\g_N=\oa_{2n}$. Then the coefficients of all polynomials
$l_{m\tss i}(z)$ and \eqref{pfexp}
generate a commutative subalgebra of $\U(\g_N)$.
If $B$ is a regular element of $\h^*_N$, then
this subalgebra is maximal commutative.
It is freely generated by the elements $l^{\tss(1)}_{2k\ts 2k},
\dots,l^{\tss(2k)}_{2k\ts 2k}$ with $k=1,\dots,n-1$
and $p^{(1)},\dots,p^{(n)}$.
\qed
\eco

\medskip

Consider the rational $R$-matrix
$R(z)=1-Pz^{-1}+Q(z-\vk)^{-1}$
associated with the Lie algebra $\g_N$, where $\vk$ equals
$N/2-1$ or $N/2+1$ in the orthogonal and symplectic case,
respectively, and $P$ and $Q$ are defined in \eqref{pmatrix}
and \eqref{rprimeq}; see \cite{zz:rf}.
The {\it extended Yangian\/}
$\X(\g_N)$ is defined
as an associative algebra with generators
$t_{ij}^{(r)}$, where $1\leqslant i,j\leqslant N$ and $r=1,2,\dots$,
satisfying certain quadratic relations;
see \cite{d:qg}, \cite{rtf:ql}.
To write them down, set
\ben
T(z)=\sum_{i,j=1}^N e_{ij}\ot t_{ij}(z),\qquad
t_{ij}(z)=\de_{ij}+\sum_{r=1}^{\infty}t_{ij}^{(r)}\ts z^{-r}
\in\X(\g_N)[[z^{-1}]].
\een
The defining relations for the algebra $\X(\g_N)$ are
written as the $RTT$ relation
\beql{RTT}
R(z-v)\ts T_1(z)\ts T_2(v)=T_2(v)\ts T_1(z)\ts R(z-v),
\eeq
with the subscripts indicating the corresponding
copies of $\End\CC^N$ in the tensor product algebra
$\End\CC^N\ot\End\CC^N\ot \X(\g_N)$.
The {\it Yangian\/} $\Y(\g_N)$ is defined as the quotient
of the extended Yangian $\X(\g_N)$
by the relation
$T^{\tss\prime}(z+\vk)\ts T(z)=1$; see \cite{aacfr:rp},
\cite{d:qg}.

Fix an arbitrary $N\times N$ matrix $C$ over $\CC$
with the property $C\tss C'=1$.
Consider the tensor product
\ben
\End\CC^N\ot\dots\ot\End\CC^N\ot\Y(\g_N)[[z^{-1}]]
\een
with $m$ copies of $\End\CC^N$ and introduce the formal
series $\tau_m(z,C)$ with coefficients in the Yangian by
\ben
\tau_m(z,C)=\tr\ts S^{(m)}C_1\dots C_m\ts
T_1(z)\dots T_m(z+m-1)
\een
in the orthogonal case, and by
\ben
\tau_m(z,C)=\tr\ts S^{(m)}C_1\dots C_m\ts
T_1(z)\dots T_m(z-m+1)
\een
in the symplectic case with $N=2n$ and $m\leqslant n$.

\bpr\label{prop:bs}
The coefficients of the series $\tau_m(z,C)$, $m\geqslant 1$,
generate a commutative subalgebra of $\Y(\g_N)$.
\epr

\bpf
The argument is based on the fact that the symmetrizer
$S^{(m)}$ admits a multiplicative presentation \eqref{symfpnew}
and uses standard $R$-matrix techniques;
see e.g. \cite[Sec.~1.14]{m:yc}. First suppose that $C=1$.
The $RTT$-relation \eqref{RTT} implies
\beql{funda}
R(z_1,\dots, z_{m+l})\, T_1 (z_1) \dots T_{m+l} (z_{m+l})
= T_{m+l} (z_{m+l}) \dots T_1(z_1)\, R(z_1,\dots, z_{m+l}),
\eeq
where the $z_i$ are formal variables and
\beql{prrm}
R(z_1,\dots, z_{m+l})=\prod_{1\leqslant i<j\leqslant m+l}
R_{ij}(z_i-z_j)
\eeq
with the product taken in the lexicographic order
on the pairs $(i,j)$. Now specialize the variables by setting
\ben
z_i=z+i-1,\quad i=1,\dots,m \Fand
z_{m+j}=v+j-1,\quad j=1,\dots,l
\een
in the orthogonal case, and by
\ben
z_i=z-i+1,\quad i=1,\dots,m \Fand
z_{m+j}=v-j+1,\quad j=1,\dots,l
\een
in the symplectic case, where $z$ and $v$ are formal variables.
The product \eqref{prrm} then becomes
$
R(z_1,\dots, z_{m+l})=\wt R(z,v)\ts S^{(m)}S^{(l)\tss\prime},
$
where $S^{(l)\tss\prime}$ is the image
of the Brauer algebra symmetrizer in the tensor
product of the copies of $\End\CC^N$ labeled by $m+1,\dots,m+l$,
and $\wt R(z,v)$ is the product of the $R$-matrix
factors of the form $R_{ij}(z_i-z_j)$ with
$1\leqslant i\leqslant m$ and $m+1\leqslant j\leqslant m+l$.
Since all these factors are invertible, using \eqref{funda}
and taking the trace over all copies
of $\End\CC^N$ we find that
\ben
\tau^{}_m(z)\ts\tau^{}_l(v)=\tr\ts\wt R(z,v)^{-1}
\tau^{}_l(v)\ts \tau^{}_m(z)\ts \wt R(z,v)=\tau^{}_l(v)\ts \tau^{}_m(z).
\een
The $RTT$ relation \eqref{RTT} will hold if
$T(z)$ is replaced with the matrix $C\tss T(z)$.
Therefore the above arguments also apply to this matrix
instead of $T(z)$.
\epf

Using \eqref{funda} with $l=0$ together with the
multiplicative formula \eqref{symfpnew},
we obtain that, up to a shift $z\mapsto z+\text{const}$,
both in the orthogonal and symplectic case,
the series $\tau_m(z,C)$ can be written as
\ben
\tr\ts S^{(m)}
C_1T_1(z)\dots C_mT_m(z-m+1)=\tr\ts S^{(m)}
C_1T_1(z)\tss e^{-\di_z}\dots \ts C_mT_m(z)
\tss e^{-\di_z}\ts e^{m\tss\di_z}.
\een
The Yangian $\Y(\g_N)$ admits a filtration defined
on the generators by $\deg t_{ij}^{(r)}=r-1$. The associated
graded algebra $\gr\Y(\g_N)$ is isomorphic to the universal
enveloping algebra $\U(\g_N[t])$; see \cite[Theorem~3.6]{amr:rp}.
Following \cite[Sec.~10]{mtv:be}, extend the filtration to
the algebra $\Y(\g_N)[[z^{-1},\di_z]]$ by setting
$\deg z^{-1}=\deg\di_z=-1$, and regard the matrix $C$ as an element
of an extended filtered algebra such
that $C-1$ has degree $\leqslant-1$
with the image $B$ in the component of degree $-1$
of the associated graded algebra. Then the matrix $B$ has the property
$B+B^{\tss\prime}=0$. The entries
of the matrix $1-e^{-\di_z}CT(z)$
have degree $-1$ and its image in the graded algebra
coincides with the matrix $\di_z-B-F(z)_-$.
Therefore, the commutative subalgebra in the Yangian $\Y(\g_N)$
provided by Proposition~\ref{prop:bs} gives rise to a commutative
subalgebra in the associated graded algebra $\U(\g_N[t])$
by the argument originated in \cite{t:qg}; see also
\cite[Sec.~10]{mtv:be} for more details.
Moreover, in the orthogonal case
this subalgebra coincides with the one
obtained in Corollary~\ref{cor:commsub}, while
in the symplectic case this subalgebra is smaller due to
the restriction $m\leqslant n$.


\end{document}